\documentclass[conference,a4paper,]{IEEEtran}

 \usepackage{mdframed} \usepackage{lipsum}

\usepackage[centertags]{amsmath}
\usepackage{amsfonts,amssymb,amsmath}
\usepackage{graphicx}
\usepackage{multirow}
\newtheorem{thm}{\bf Theorem}

\newtheorem{lemma}{\bf Lemma  }

\newcommand{\ForITCTR}{\ifthenelse{1<2}}   

\usepackage{color}

 \usepackage{ifthen}
 
  \newcommand{\up}{^}
\newcommand{\Draft}{\ifthenelse{1<0}}
\newcommand{\Dtls}{\ifthenelse{1<2}}

\newcommand{\ignore}[1]{}

\newcommand{\TR}{\ifthenelse{1>2}}

\newcommand{\NTR}{\ifthenelse{1>0}}

\newcommand{\Thesis}{\ifthenelse{1<2}}

\newcommand{\Conf}{\ifthenelse{1<0}}

 \renewcommand{\i} {{ \epsilon }}
 \renewcommand{\t}{{\tau}}

\newcommand{\Expp}{\mathbb{E}}

 \newcommand{\Id}{{\cal I}}

 \newcommand{\csn}{\check}
 
 \newcommand{\eop}{{\hfill $\blacksquare$} }
 
 \Conf{
\onecolumn
 }
 {}
 
\begin{document}
\NTR{
\title {Queuing  with Heterogeneous Users: Block Probability and  Sojourn times}}
{
\title {Static achievable region  with impatient customers}
}

 \author{Veeraruna Kavitha and Raman Kumar Sinha \\ IEOR, Indian Institute of Technology Bombay, Mumbai, India}
   
\maketitle

\newpage 
\begin{abstract}
 
Communication networks need to support voice and data calls simultaneously.  This results in a  queueing system with   heterogeneous agents.  One class of agents demand immediate service,  would leave the system if not provided. The second class of customers have longer job requirements and can wait for their turn. We discuss the achievable region of such a two class system, which is the region of all possible pairs of performance metrics. 
 Blocking probability is the relevant performance for eager/impatient class while  the expected sojourn time is appropriate for the second  tolerant class.  We obtain the achievable region, considering static  policies that do not depend upon the state of the second class. 
We conjecture a pseudo conservation law, in a fluid limit for eager customers, which relates the  blocking probability of eager customers with the expected sojourn time of the tolerant customers. Using this conjecture we obtain the static achievable region.  We validate the pseudo conservation law using two example families of static schedulers, both of which achieve all the points on the achievable region. 
  Along the way we obtain smooth control  (sharing)
  of resources between    voice and data calls.   We also consider an example dynamic policy to establish that the dynamic achievable region
  is strictly bigger than the static region, for this  heterogeneous  queueing system.
 
\end{abstract}
 
 \NTR{}{
 \vspace{-2.5mm}}
 
 \section{Introduction}
We consider queueing systems with heterogeneous classes of users. One is a class of eager/impatient customers, who would reject the system if service is not offered  immediately. Alternatively they may be willing to wait for a very brief period, but they would like to spend minimal time with system  as they have short job requirements. The other class of the customers are tolerant, can wait for their turn. One would require parallel service-offer facility to handle the first class. These customers are satisfied as long as their service starts, even if they have to share the service facility with others. However, there may be a limit to which the eager class of customers are ready to share the service utility.  
Our aim in this paper is to study the achievable region,  basically all possible pairs of `heterogeneous' performance measures of the two  classes, under certain conditions.

One of the main motivations for  this  paper is data-voice calls of a communication network. 
Data calls are delay tolerant, but require precision.  They can tolerate delays in service, but not the errors in transmission. Their job requirements are usually long.
On the other hand the voice calls are impatient, need immediate service.
However their job requirements would usually be smaller.
We consider two policies for capacity/resource sharing between the data-voice calls. In the first policy  entire capacity is transferred to the voice calls (when admitted), irrespective of the number   receiving the service.  The voice calls operate in processor sharing mode, and we refer this as $PS$ policy.   In  the second ({\large  c}apacity {\large d}ivision/$CD$) scheduling policy,   the capacity   transferred  (resources allocated) to   voice calls  is proportional to the number receiving the service. In this paper,  we  study both the policies. 

Consider a communication system with  $K$ orthogonal channels. For example, each channel could be one or a collection of resource blocks as in an OFDM based LTE network (e.g.,\cite{LTE}). 
Initially all the channels  are dedicated to data calls.  As and when the voice calls arrive, one by one the channels are transferred and data calls use the remaining.   Our $CD$ policy captures  this scenario precisely.  We provide (admission control based)  policies for such scenarios:  when admission of voice calls is severely reduced the expected time spent by data calls in the system would be less and vice versa. We show that the admission control based policies achieve the entire span of the 'achievable region' of ordered pairs of block probability and expected sojourn times. 
If   the voice calls are served at the highest possible rate as with $PS$ policy, it  improves the chances of a free server being available to subsequent voice arrivals.  The two achievable regions overlap, but  $PS$ has a bigger region (when  number of maximum parallel calls is fixed) as it    attains a   smaller blocking probability.

 \Conf{The performance metric important for impatient customers is the blocking probability, the probability that a customer returns without service. The tolerant customers can wait for their turn, however their satisfaction depends upon the expected sojourn time. This is the total time spent in the system.We are interested in obtaining the achievable region of such queueing systems.

 An achievable region for a system with $n$-classes is the set of all  possible relevant performance vectors $(pm_1, \cdots, pm_n)$, obtained by varying all possible scheduling policies (\cite{shanthikumar1992multiclass,federgruen,bertsimas} etc.,). In our heterogeneous setting, achievable region is the set of all possible pairs of blocking probability and expected sojourn times. The main focus of the paper is this achievable region.
  }
{
An achievable region for a system with $n$-classes is the set of all   possible relevant performance vectors 
$(pm_1, \cdots, pm_n)$, obtained by varying all possible scheduling policies (\cite{federgruen,bertsimas,shanthikumar1992multiclass} etc.,).
The performance metric important for impatient customers is the blocking probability, the probability that a customer returns without service. The tolerant customers can wait for their turn, however their satisfaction depends upon the expected sojourn time. This is the total time spent in the system.  This paper focuses on heterogeneous achievable region, the set of all  possible pairs of blocking probability and expected sojourn time.  Once the achievable region is known, many relevant optimization problems can be solved easily. For example, the problem of finding the optimal expected sojourn time of data calls, given a constraint on blocking probability of voice calls can readily be solved.

}
 The achievable region is well understood for homogeneous  classes, when  the performance metric of both the classes  is expected sojourn/waiting time.  
 Conservation laws, pioneered by \cite{Kleinrock1965}, capture the fundamental limits of the performance measures like mean waiting time of various classes of customers when they share a common server.  \Conf{Coffman and Mitrani \cite{coffman} were the first to identify such achievable regions when they identified it for multi-class $M/M/1$ queue with pre-emptive discipline among policies with regenerative structure.}{}
Multi-class single server   queueing systems pose nice geometric structure (polytopes) for achievable region (e.g., \cite{coffman,shanthikumar1992multiclass}).
\Conf{
A Multi-Class,  Multi-Server Queue with preemptive priorities (i.e. a high priority customer preempts the lower one,if all the servers are serving low priority customers) of two groups, each group may consist of multiple customer types, has been discussed in (\cite{Sleptchenko}). However in (\cite{Harchol}) authors introduce Recursive Dimensionality Reduction (RDR) technique for multi-server queueing systems with multiple priority classes. An aproximation of sojourn time is derived in (\cite{Izagirre}) for multi-class single-server queue with $K$ classes of customers, where the service capacity is shared simultaneously among all customers present in proportion to the respective class-dependent weights.


Another community of researchers focused on dynamic control over multi-class queueing systems due to its various applications in computers, communication networks, and manufacturing systems. One of the main tools for such control problems is to characterize the achievable region, then use optimization methods to find optimal control policy (see \cite{bertsimas1995achievable}, \cite{bertsimas1996conservation}, \cite{li2012delay} etc.,).  

A parametrized family of scheduling policies is called \textit{complete} by \cite{complete}, if it achieves all possible performance vectors of the achievable region, average waiting times in their context.   A complete scheduling class can be used to find the optimal control policy over all scheduling disciplines. Discriminatory processor sharing (DPS) class of parametrized dynamic priority schedulers is identified as a \textit{complete} family in case of two class $M/G/1$ queue  in \cite{hassin2009use}.   Many more families of scheduling policies are identified to be complete.

    Models with heterogeneous classes can have extensive applications in  communication networks. 
    For example voice calls in a cellular network are the impatient customers, while data calls are  the tolerant customers.
      Voice calls will be dropped in case of server unavailability, while the data calls can wait. Primary and secondary users in Cognitive networks could be another  interesting   example.   We further discuss this topic   after the model is completely introduced. 
      
}{

We conjecture a relation between the expected sojourn time and the blocking probability  (for any  static scheduling policy), and, \underline{call it a pseudo conservation   law}.   This pseudo conservation law is valid in a fluid limit for short job impatient agents.  We then show that  two sets of scheduling ($PS$ and $CD$) policies satisfy this conservation law and also achieve all the points of the resulting achievable region. 
   
   To the best of our knowledge we are not aware of a work that directly studies this type of  a heterogeneous achievable region. Some variants of queueing systems  (e.g., \cite{Whitt,Harchol,Izagirre,CD_close,Roland,Sleptchenko,White} etc) have some connections to few parts of  our models, and these are discussed   in \cite{Journ1}.  
} 

In \cite{CD_close} authors consider multi-class queueing system with eager and tolerant customers.  This  is the queueing system  that is closest to the one considered in this paper, especially to  $CD$ policy. 
With our CD policy, the tolerant customers utilize all the remaining servers,   
and hence the system is  work conserving (only) with regard to the tolerant customers.
While in  \cite{CD_close}   tolerant customers are also served in multi-server mode, i.e.,  each tolerant customer is provided one server and the idle servers (if any) are not utilized.
Further the authors in  \cite{CD_close} obtain a set of (balance) equations, solving which  stationary probabilities (and then the stationary  performance) can be derived. We provide a closed form expression for these performance measures in fluid limit for eager customers.  We also conjecture  a `policy-independent' pseudo conservation law.  Generally it is easy to obtain the performance of the loss systems (customers are lost due to their impatience), as these systems are usually finite state space-Markov chains.  However the challenging part is to obtain the performance for the tolerant customers. One can obtain system of equations (as in \cite{CD_close}), the solution of which provides the performance measure. However once the pseudo conservation law is proved, the performance of tolerant customers is immediately known as soon as the performance of the eager customers is computed.

There has been considerable work that discusses resource sharing between voice and data calls and we discuss a few here. 
In \cite{Tang} authors consider a three channel pool scheme, and obtain a novel adjustable boundary
based channel allocation scheme with pre-emptive
priority for integrated voice and data  networks.   They attain various levels of priority by adjusting the division  of the total available channels among the three pools. 
In \cite {Zhang} authors again consider channel allocation scheme for packet level allocations.  These papers {\it  discuss coarse sharing of resources between data-voice calls.}
In our model \Conf{we consider continuous (exponential) job requirements and }{we }provide a scheme to smoothly control the performance measures of the two classes. By varying the admission parameter,  smoothly in the interval $[0,1]$, one can achieve any pair of performance measures on the static achievable region.

We predominantly discuss  (partially) 'static policies', the scheduling policies  that does not depend upon the state of the tolerant customers. Towards the end, we consider an example dynamic policy and demonstrate that the achievable region with dynamic policies is strictly bigger.

Typically 
customers with long job requirements form tolerant class, while the  eager ones  demand for short jobs.  Our results are applicable in an asymptotic fluid regime that takes advantage of this.
 In section \ref{sec_model} we describe the problem and pseudo conservation law is conjectured in section \ref{sec_conjecture}.  $PS$ and $CD$ policies are respectively discussed in sections \ref{sec_static} and  \ref{sec_CD}.   
 And dynamic policies are discussed in section \ref{sec_dynamic}.
  
\section{Problem Statement and System model}
\label{sec_model}
 \vspace{-1mm}
 \ignore{We consider a heterogeneous  queueing system with  two classes  of customers and possible parallel service capabilities. 
 The first class of customers  have short job requirements  and are impatient  in nature. They either demand immediate service or service after a short wait period.  More realistically they would  like to spend minimal time with the system. While the second class of customers have considerably long  job requirements and are wiling to wait for their turn. }
 We refer the impatient/{\large  \bf e}ager customers by $\i$-customers while the {\large \bf  t}olerant customers are referred to as $\t$-customers.
The system has a fixed server capacity, that needs to be shared between the two classes of customers.  
The exact sharing of capacity depends upon the allocation/scheduling policy. For example the system can serve $K$ customers in parallel for some $K$, by  dividing the  server capacity among the customers under service.  The system can chose to vary $K$ dynamically, e.g., processor sharing. The system can chose to serve one customer with full capacity etc.   In this paper we discuss two example sets of scheduling policies.  
{\it We only consider  `$\t$-work conserving' policies, wherein the $\t$-customers utilize all the remaining server capacity.}

\noindent{\it Arrival process and the Jobs:} 
 The arrival processes are modelled by independent Poisson processes,   with rates $ \lambda_{\i}$ and  $\lambda_{\t}$ respectively.    
 The job requirements for both the classes are   exponentially distributed.  The time required to complete a job, depends upon the scheduling policy. If  a $\t$-customer ($\i$-customer) is served with full server capacity, then the service time  is exponentially distributed with parameter   $\mu_\t$ (respectively   $\mu_\i$).


 \vspace{-1mm}
 
 \subsection{Achievable region}
 
 The two classes of users have different goals and hence naturally   require  different qualities of service (QoS). 
An eager $\i$-agent would leave the facility without service, if service is not provided almost immediately.
  Hence a scheduling decision  (basically admission decision) is required at every   $\i$-arrival instance\footnote{On the contrary, in homogeneous setting two or more classes of agents wait at their waiting lines and scheduling epochs are the service completion/departure epochs. The scheduler had to decide which class  to be  served next. While in heterogeneous setting,  at any departure epoch there is only one class of agents possibly waiting and hence no decision is required.}. 
   {\it Block probability $P_{B}  $, the probability of customers departing without  service,}  is an important performance metric for $\i$-class. 
This also implies the service of a typical $\t$-customer can possibly be interrupted, possibly to provide required QoS for $\i$-agents,  and a typical $\t$-agent can face several such interruptions during its service.  Thus  the {\it expected sojourn time $E[S_\t ]$, the expected value of the total time spent by a typical agent would be an appropriate QoS for $\t$-class.}  It is not sufficient to consider the expected waiting time, the time before the service start of a typical $\t$-agent. 
   Either of these performance metrics depend upon the scheduler $\beta$ used. 
   Therefore the achievable region is given by: 
   
   \vspace{-2mm}
{\small  $$ {\cal A}_{hetero}   = \lbrace \left( P_{B}(\beta),E[S_\t(\beta)] \right) : { \beta  } \mbox{ is a scheduler}  \rbrace .$$}
 In this paper {\it we consider  ($\t$) static policies, wherein the  $\i$- admission rules do not depend upon  the status  of the $\t$-class.}   The probability of admission, $p$,    is an important parameter of any such scheduling policy.
Further,   the (maximum) number of $\i$-calls  served in parallel and the  sharing of resources between $\i$ and $\t$ customers is also a part of the scheduling decision.     For example, the system may allocate/transfer entire server capacity to the first admitted $\i$-arrival. It may processor-share the capacity among the further admitted $\i$-customers. There may be a limit on the number of $\i$-customers that can simultaneously share the capacity. Alternatively 
the system may allocate 
 a fixed fraction  of the server capacity to each admitted $\i$-arrival and the remaining is allocated to $\t$-class etc. All these rules are independent of the   $\t$-state (e.g., the number of $\t$ customers in the system, waiting time of them etc).   This implies that $\i$-calls pre-empt  $\t$-call when required.  
 In all  a $\t$-static policy implies that  an $\i$-arrival is admitted with some probability $p$, and further  admission also depends upon the number of $\i$-calls already in the system, but not on $\t$-state. Mathematically a static achievable region is defined:  

\vspace{-4mm}
{\small \begin{eqnarray}
\label{eqn_hetero_achieve}
  {\cal A}_{hetero}^{static} =  \Big \lbrace \left( P_{B}(\beta^{CS}_p),E[S_\t(\beta^{CS}_p)] \right) \mbox{\large  : }\hspace{24mm} \nonumber \\  0 \le p  \le 1 \mbox{  and  $(CS)$ a capacity sharing rule }  \Big \rbrace  .
 \end{eqnarray}}
We primarily analyze the static achievable region.  In \NTR{section \ref{sec_dynamic}}{\cite{Journ}} an example  dynamic policy (also depends upon $\t$-state)  is considered to show that the achievable region with dynamic policies is strictly bigger than the static achievable region.

 \vspace{-2mm}

\subsection{Short-Frequent Job (SFJ) limits}  
The $\i$-class  has short job requirements. If one considers limit $\mu_\i \to \infty$, the impact of $\i$-customers becomes negligible  at the limit. 
To obtain a more general and useful result, we also increase the $\i$-arrival rate   while $\mu_\i \to \infty$.  
That is, every $\i$-agent may utilize  the server  for a short duration, but the system has to attend  the $\i$-agents  frequently. 
Because of this $\i$-agents cause  significant impact even in  the limit.  
To be more precise we consider the limits $\mu_\i \rightarrow \infty\; \mbox{and}\; \lambda_\i \rightarrow \infty$ while the load factor $\rho_\i =\lambda_\i/\mu_\i $ is maintained constant. We  refer this  as  ``Short-Frequent Job (SFJ) limits''.

\ignore{
\subsection*{Remarks about the model}

One of the motivations for  this model is   data-voice calls of communication networks. 
Data calls are delay tolerant, however require precision.  
 They prefer  transmission using   highest possible resources to ensure quality.  Further often they require transfer of long files. Thus the job requirements of a typical data call is usually long.    On the other hand the voice calls are impatient, however the ear can tolerate inaccuracies to certain extent. Hence the data corresponding to  voice users is sufficiently compressed, resulting in short job requirements. 
}

\section{Conjecture of a Pseudo Conservation law}
\label{sec_conjecture} 

In a   multi-class queueing system  with all tolerant classes  (homogeneous system),  a work conservation law holds. 
 The total workload in the system remains the same irrespective of the scheduling policy, as long as  the server does not idle  during busy period.   Further by Little's law and Wald's lemma, a linear combination of expected sojourn (or waiting) times of different classes of customers remains the same irrespective of the scheduling policy (e.g., \cite{shanthikumar1992multiclass}).  

The above is obviously true when   the incoming workload remains the same.  
However in our heterogeneous setting, the $\i$-customers depart the system, if service is not offered immediately.  And this depends upon the scheduling policy. Thus the workload arriving into the system itself changes with different scheduling policies and naturally one may not expect work conservation.  However  if the amount of work blocked remains the same, one can   anticipate a different kind of work conservation. 
We conjecture that  given a probability of blocking, irrespective of the way the $\i$-agents are blocked and irrespective of the way the $\t$-agents are served, the $\t$- expected sojourn time   remains the same\footnote{Within tolerant class of customers, the expected sojourn time by Little's law and Wald's lemma is proportional to the  workload in the system}.   And this could be conjectured only in SFJ limit and when the policies do not depend upon the $\t$-state (The proof of this conjecture is considered in \cite{Journ1}).  

 In SFJ limit, $\i$-agents   will have fluid arrivals and departures.  
 Given  the $\i$-load factor ($\rho_\i$) and 
 the probability of blocking  ($p_B$),   in the  SFJ limit, the $\i$-agents occupy   $\rho_\i (1-p_B)$  fraction of  system resources at all the times.  Hence we 
  {\it conjecture } that   the $\t$- performance equals that of an  $M/M/1$ queue with smaller service rate $\mu_\t (1-\rho_\i (1-p_B))$, and  that   the expected sojourn time for any $0 \le p_B \le 1$  equals:
 
 \vspace{1mm}
 \hspace{-2mm}
\hspace{-2mm}\fbox{\begin{minipage}{8.7cm}\vspace{-2.5mm}
{\Draft{}{\small}\begin{equation}
 E_{S_\t} (p_B) :=   \frac{1}{\mu_\t (1-\rho_\i (1-p_B)) - \lambda_\t} \mbox{ if }     \rho_\i (1-p_B) + \rho_\t  < 1.  \label{Eqn_PC}
\end{equation}}\vspace{-4mm}
\end{minipage}
}

\vspace{1.5mm}
{\bf Conjecture:}  Static achievable region, in SFJ limit, equals:
 {\Draft{}{\small} \[
\boxed{ {\cal A}^{hetero}_{static} = \Big \{  \big (p_B, E_{S_\t} (p_B)  \big ) \, \Big | \,  {p_B \in [0,1], \, \,} \rho_\i (1-p_B) + \rho_\t  < 1\Big \}.}
 \]}
 We would like to refer the 
  equation (\ref{Eqn_PC}) as a {\it pseudo conservation law,} as it provides the expected sojourn time in terms of the fraction blocked (lost).  
 This would require an explicit proof which is considered in \cite{Journ1}.  
For now,   we consider two example families of schedulers and illustrate the validity of our conjecture. 
Further,   using the same sets of schedulers, we  achieve all the points of  the static region.  Such a family is generally {\it referred to as complete family of schedulers.  }

\section{Processor sharing  $PS-(p, K)$ schedulers} 
\label{sec_static}
Any $\i$-arrival is admitted to the system with probability $p$, independent of  $\t$-state. Once admitted  it will pre-empt the existing $\t$-agent, if any.   
 We consider $K$-processor sharing service discipline for   $\i$-agents. If there is only one agent of the  $\i$-class receiving service, it is served with maximum capacity, i.e., using   capacity  ${\mu}_\i.$  Upon a new (admitted) arrival of the same class, the capacity is shared among the two.   Both are served in parallel and independently, each with rate ${  \mu}_\i/2$. Upon a third (admitted)  arrival  each is served with rate ${  \mu}_\i/3.$ This continues up to $K$ $\i$-agents.  Any further   arrival,   leaves without service even after being admitted.  When any of the existing $\i$-agents depart, the service rate is readjusted to an appropriate higher value.  The $\t$-service is resumed only after all the $\i$-agents depart. We call this as $\beta^{PS}_{p, K}$ scheduling policy.
 
Tolerant  agents are served in FCFS (first come first serve) basis. They are   served in a serial fashion and with full capacity\footnote{Capacity of the server is such that, it can either serve one tolerant agent at rate $\mu_\t$, or   $l$ 
$\i$-agents each at $\mu_\i/l$ (where $l \le K$).} $\mu_{\t}$.  That is, system would serve at maximum one $\t$-agent, and the 
service of the next $\t$-agent begins only after the  preceding one departs.

 The   transitions and evolution of the $\i$-agents is independent of that of $\t$-agents    under a  static  policy:  the arrivals are admitted  and the service  is  provided to the admitted agents immediately, 
 irrespective of the state of $\t$-agents.  Thus one can analyze the $\i$-class independently and we first consider this analysis. 
 \vspace{-3mm}

\subsection{Blocking Probability of $\i$-class} \label{sec_Pb}
Fix  $0 \le p \le 1$, $K$ and consider policy  $\beta_{p,K}^{PS}$. Blocking probability is the probability with which a new ($\i$-class) arrival   leaves the system without service.   Blocking can occur in case of two events. Upon arrival, an $\i$-agent is admitted to the system with probability $p$ and is    blocked with  probability $(1-p)$.  Secondly, an admitted agent leaves   without service, if  the system is already serving $K$  $\i$-agents.

Let $\Phi_\i (t)$ represent  the number of $\i$-agents in the system  at time $t$.
 We claim that   the $\i$-class transitions are caused by  exponential random events and hence that  $\Phi_\i (t)$  is a continuous time Markov jump process (see for example \cite{Hoel}) for the following reasons:  a) it is clear that the inter-arrival times are exponentially distributed with parameter $\lambda_\i p$; b)   by Lemma \ref{lemma_exp}, given below,  the departure times are exponentially distributed with parameter $\mu_\i$  (i.e., $\sim \exp(\mu_\i$)), irrespective of  state $\Phi_\i (t)$. 

\begin{lemma}
\label{lemma_exp}
    Let $D\up{l}_\i$ represent the time to first departure among the $l$ $\i$-agents receiving the service, with {\small$1 \le l \le K$}. Then for $PS$ policy, $D\up{l}_\i\sim \exp(\mu_\i)$  for any $l$.
    \end{lemma} 
{\bf Proof:}  When $l$ agents are receiving service in parallel,  because of processor sharing  the service time of each   is exponentially distributed with parameter $\mu_\i/l$.  And the time to first departure, the minimum of these $l$ exponential random variables,  is again exponential with parameter $l  \mu_\i/l = \mu_\i$.  \eop

 \begin{figure}[h]
\vspace{-19mm}
\hspace{-7mm}
\includegraphics[width=11.3cm,height=12cm]{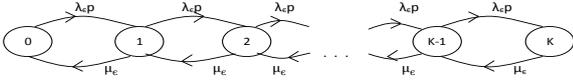}  
\vspace{-97mm}
\caption{ State transitions  for $\i$-agents with $\beta^{PS}_{p, K}$ scheduler. \label{Figure_Transitions}}
\vspace{-2mm}
\end{figure}

  In Figure  \ref{Figure_Transitions}, we depict the transitions of the continuous time Markov jump process $\Phi_\i (t)$. For such processes, well known balance equations are solved to obtain the stationary probabilities (see for example \cite{Hoel}). The stationary probabilities, $\{\pi_0, \pi_1, \cdots, \pi_K\}$, of  $\Phi_\i (t)$ are obtained by solving:  

\vspace{-4mm}
 {\Draft{}{\small}\begin{eqnarray*}
  \pi_0 \lambda_\i p\Draft{&=&}{\hspace{-2mm}  &=&\hspace{-2mm} }  \mu_\i \pi_1, \Draft{\hspace{5mm}}{\hspace{1.2mm}} 
\pi_l ( \lambda_\i  p + \mu_\i)  = \lambda_\i p \pi_{l-1} + \mu_\i \pi_{l+1}   
   \mbox{ for  }   1\le l < K, \\
 \mbox{ and } &&
    \pi_K  \mu_\i  =  \lambda_\i p \pi_{K-1}.  
  \end{eqnarray*}}
 The solution or the stationary probabilities are ($0\le l \le K $): 

\vspace{-5mm}
 {\Draft{}{\small} \begin{eqnarray*}  
\pi_l = \frac{  \rho_{\i,p}^l }{a_0 }  
 \,\, \mbox{ with }   
  \,\, a_0 :=  \sum\limits_{j=0}^{K}\rho_{\i,p}^{j} \mbox{ and }
\rho_{\i,p} := \frac{\lambda_{\i}p}{\mu_{\i}} = \rho_\i p. 
\end{eqnarray*}}
  An admitted agent gets blocked, if it finds $K$ $\i$-agents in the system, and,  this  by PASTA (Poisson Arrivals See Time Averages)  equals the stationary probability $\pi_K $ of $K$ $\i$-agents in the system.  The  agents  are not admitted   with probability $(1-p)$ and those admitted are blocked with probability $\pi_K$. Therefore the overall blocking probability equals: 
  
  \vspace{-4mm}
 {\Draft{}{\small} \begin{eqnarray}
\boxed{ P_{B}^{PS}(p) = (1-p) + p \pi_K 
= (1-p)  +  p\frac{\rho_{\i,p}^{K}}{a_0} .}
\label{Eqn_PBp}
\end{eqnarray}}
\vspace{-6mm}
\subsection{Expected sojourn time of $\t$-class} 
\label{sec_E_t}
The $\i$-class requires short but frequent jobs (e.g.,  voice calls).  Hence we are looking for a good relevant approximation  that  facilitates the analysis, and which further allows us to study other important variants (like  $CD$ policy of section \ref{sec_CD}). 
Towards this, we approximately (accurate  asymptotically) decouple the evolution of $\t$-agents from that of $\i$-agents.    

We first understand the effective server time (EST),  $\Upsilon_\t$, which is defined as the total time period between the service start and the service end of a typical $\t$-agent. 
We refer this as EST of the agent under consideration, as no other $\t$-agent has access to server during this period. 
Sojourn time of a typical $\t$-agent equals the sum of two terms:  a)  waiting time, the time before the service start; and  b) EST $\Upsilon_\t$, the time after the service start. 
 
\subsubsection{Analysis of effective server time (EST) $(\Upsilon_\t)$}
 This time equals the sum of the actual service time, $B_\t$, of the $\t$-agent and the overall time of interruptions  caused by $\i$-agents, which is denoted by  $\Upsilon^e_\t$.  Let  $N(B_\t)$ represent the total number of  the $\i$-class interruptions, that occurred during the service time $B_\t$. 
  In reality  these interruptions would have occurred in disjoint time intervals, the sum of all of which is $B_\t$. This random number has   same stochastic nature as the number of Poisson arrivals that would have occurred in a continuous time interval of length $B_\t$. This is true because of the memory less property of the exponential service time  $B_\t$ and because Poisson process is a counting process. After an $\i$-agent interrupts the  ongoing  $\t$-agent, there is a possibility of further admissions. Eventually the service of the $\t$-class is resumed, where  left, when all the $\i$-agents (that were admitted) leave the system.  

\begin{figure*}
\vspace{2mm}
\begin{minipage}{11cm}
\vspace{-13.8mm}
\hspace{-9mm}
\includegraphics[width=11cm,height=11cm]{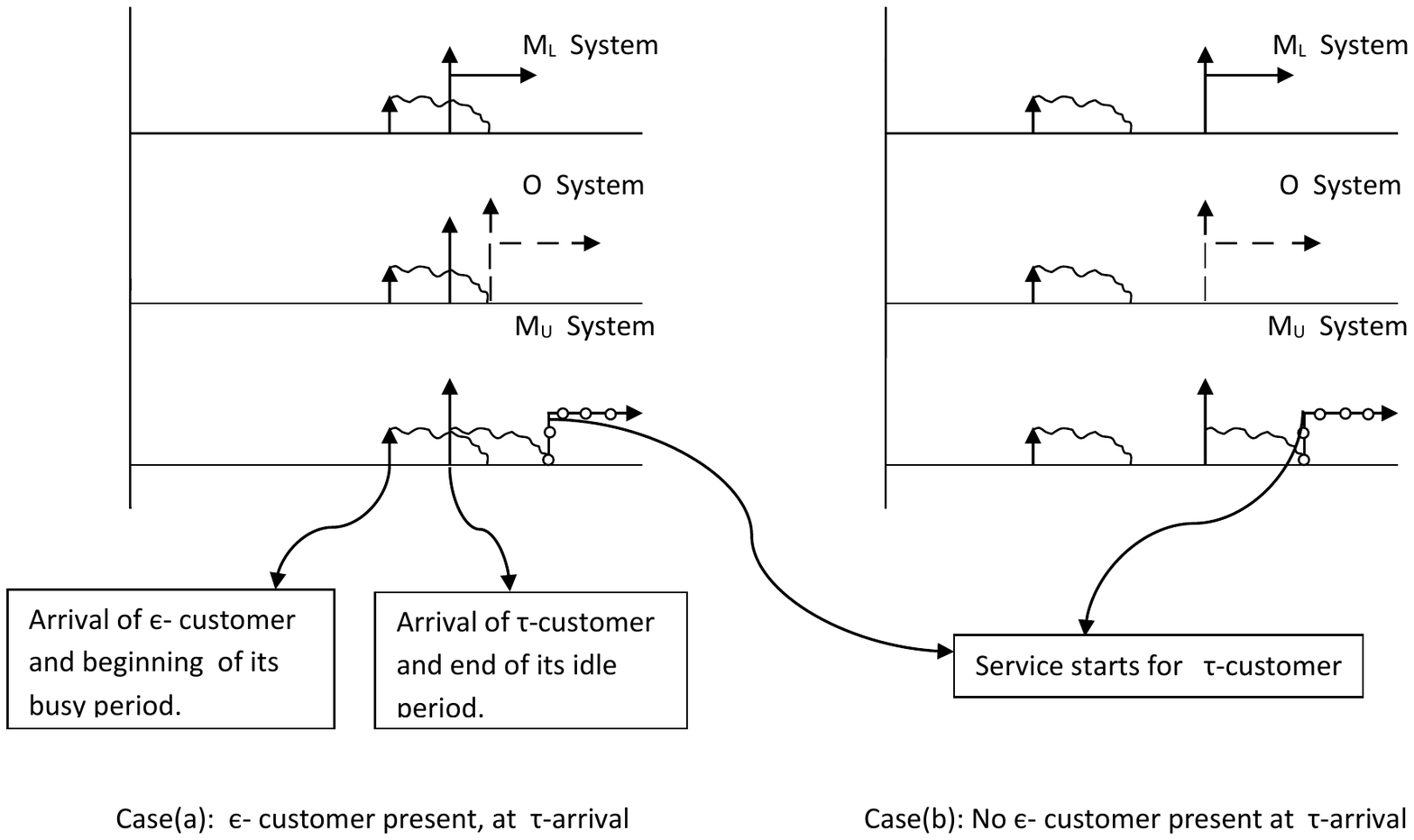}  
\vspace{-60mm}
\caption{ Example sample paths of the three systems ($PS$ policy) \label{Figure_Threesystems}}
\end{minipage}
\hspace{-22mm}
\begin{minipage}{10cm}
 
\vspace{-37mm}
\hspace{-2mm}
\includegraphics[width=8cm,height=9.cm]{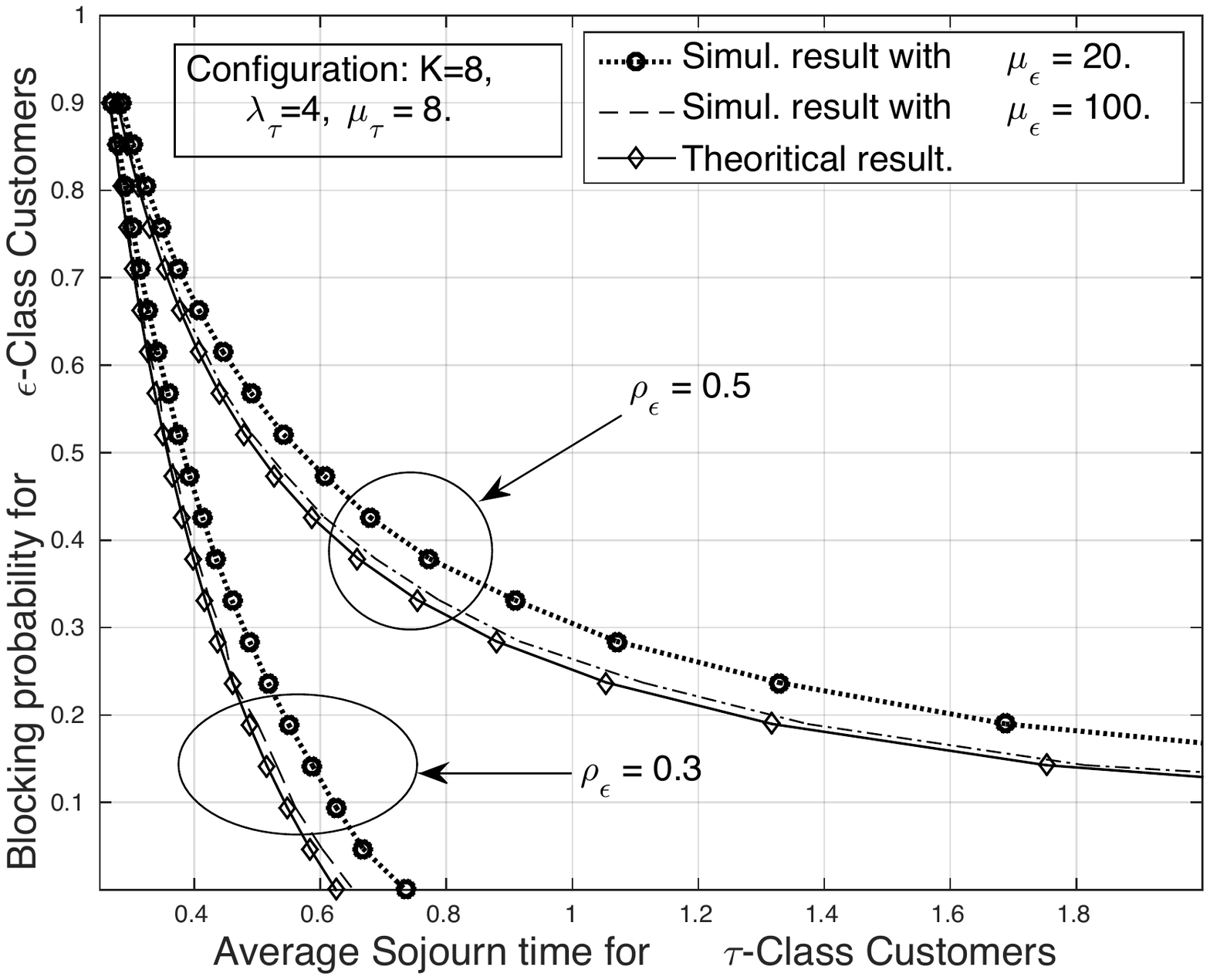} 
 
\vspace{-27mm}
\caption{  Achievable region:  Simulated and Theoretical results \label{Figure_Sim_thARegion}}
 
\vspace{-5mm}

\end{minipage}

\vspace{-5mm}
\end{figure*}
 
  Thus the time duration for which the service of $\t$-agent is suspended per interruption, equals a busy period of the $\i$-class,   that started with one $\i$-agent. There would be $N(B_\t)$ (random) number of such interruptions.  Hence,
  
  \vspace{-5mm}
{\Draft{}{\small}\begin{eqnarray}
\Upsilon_{\t} = B_{\t} + \Upsilon^{e}_{\t} \mbox{ with } \Upsilon_{\t}^e := \sum\limits_{i=1}^{N(B_{\t})} \Psi_{\i,i} \, ,\label{Eqn_Tt}
\end{eqnarray} }where
  $\{\Psi_{\i,i} \}_i$ are the  IID (independent and identically distributed) copies of $\i$-busy period.  
  We have the following result. 
  \begin{lemma}
\label{Lemma_Tt}
  The first two moments of the $\i$-busy period and   EST $\Upsilon_\t$ are given by:   
  
  \vspace{-5mm}
{\Draft{}{\small}\begin{eqnarray} \label{Eqn_EPsiUpsilon}
E[\Psi_\i] \hspace{-3mm}&= &\hspace{-3mm} \frac{a_1  }{\mu_{\i}} \mbox{ and }   
 E[\Psi_\i^2 ] =  \frac{1}{\mu_{\i}^2}  \sum_{i=1}^{K} \frac{q^{i-1}\:c_{i}}{(1-q)^{i}}  , \\
  E[\Upsilon_{\t}] \hspace{-3mm}&= &\hspace{-3mm} \frac{1}{\mu_{\t}} + \frac{ \lambda_{\i}p} {\mu_{\t}} E[\Psi_\i]
\Draft{=    \frac{1}{\mu_\t} + \frac{ \rho_{\i,p} a_1 } {\mu_\t}
}{}
   \, = \, \frac{ a_0} {\mu_\t }, \nonumber  \\
E[\Upsilon^{2}_\t]
%
\Draft{\hspace{-3mm}&= &\hspace{-3mm} \frac{2}{\mu_\t^2} \bigg [ 1+ 2 \rho_{\i, p} a_1 +  \rho^2_{\i, p} a_1^2   \bigg] + \frac{\rho_{\i,p} }{\mu_\t \mu_\i} \sum_{i=1}^{K} \frac{q^{i-1}\:c_i }{(1-q)^i },  \nonumber \\ }{}
\hspace{-3mm}&= &\hspace{-3mm} \frac{2 a_0\up2}{\mu_{\t}^2} + \nonumber \frac{\rho_{\i,p} }{\mu_\t \mu_\i} \sum_{i=1}^{K} \frac{q^{i-1}\:c_{i}}{(1-q)^{i}},
\hspace{-2mm} 
\end{eqnarray}}where the constants $q$, $\{c_i\}$ and $\{a_i\}$ are defined as:
{\Draft{}{\small}
\begin{eqnarray}
\rho_{\i,p} \hspace{-3mm}&= &\hspace{-3mm} \frac{\lambda_\i p }{ \mu_\i},  \, \, q =  \frac{\rho_{\i, p}}{ \rho_{\i,  p}+ 1}, \,\,
a_{i} =   \sum\limits_{j=0}^{K-i}\rho_{\i,p}^{j} \mbox{ for all } 0 \le i \le K, \hspace{6mm}
\label {Eqn_constants} \\
  b_{i} \hspace{-3mm}&= &\hspace{-3mm}\sum\limits_{j=K-i+1}^{K-1} (K-j)\rho^{j}_{\i,p} \mbox{ for all } 2 \le i \le K ,     \,\,\, \,  b_{1} = 0,  \nonumber \\
 c_1  
\hspace{-3mm}&= &\hspace{-3mm} \frac{2\rho_{\i,p}   \left ( 2 a_2+b_2 \right )  +2 }{(1+ \rho_{\i, p} )^2 \mu_\i^2},   \mbox{ and  for all } 1\le i < K  \nonumber \\
 c_i   
 \hspace{-3mm}&= &\hspace{-3mm} \frac{2\rho_{\i,p}   \left ( (i+1)  a_{i+1}+b_{i+1} \right )+2 \left( (i-1)a_{i-1}+b_{i-1} \right ) +2 }{(1+ \rho_{\i, p} )^2 \mu_\i^2 }, \nonumber \\
 c_K \hspace{-3mm}&= &\hspace{-3mm}  \frac{2 \rho_{\i, p}  (Ka_K+b_K) + 2 ((K-1)a_{K-1}+b_{K-1}) + 2 }{(1 + \rho_{\i,p}   )^2\mu_\i^2 } .\nonumber
 \end{eqnarray}}
\end{lemma}
\noindent
\textit{Proof:} The proof is provided in Appendix A.  {\hfill $\blacksquare$}

  \medskip
 \subsubsection{Approximate decoupling via Domination}  Every $\t$-agent undergoes similar stochastic behaviour, as below. Each  agent has to wait for the beginning of its service, and has to finish its service in the midst of random interruptions, all of which have identical stochastic nature. Further, evolution of the $\i$-agents during the EST $\Upsilon_\t$ of one $\t$-agent is  independent of that of the other $\t$-agents. Hence
 the  $\Upsilon_\t$ times corresponding to different $\t$-agents are independent of each other. 
  Thus the idea is to model the $\t$-class evolution approximately as an independent process, with that of an $M/G/1$ queue. The arrivals remain the same,  but the service times in $M/G/1$ queue are replaced by the sequence of ESTs $\{\Upsilon^t_\t\}$. 
 
 We call this $M/G/1$ queue as ${\cal M}_L$ system  and the original system as ${\cal O}$ system.
  In fact we will define another $M/G/1$ system ${\cal M}_U$ as below and show that: a) the performance  (expected sojourn times)  of the original system is bounded between the performances of the two $M/G/1$ systems; and  b) that the performances of the two sandwiching systems converge towards each other as $\mu_\i \to \infty$ (even with $\rho_\i$ fixed).

\paragraph{${\cal M}_L$ system}
The ESTs are considered as service times of $\t$-agent   in ${\cal M}_L$ system.
We study the (sample path wise) time evolution of the two  systems, original and ${\cal M}_L$, to  demonstrate the required domination. 
Towards this, we assume that both the systems are driven by same input (arrival times and service requirements) processes.  
Consider that   both the systems start with  same number   (greater than 0) of $\t$-agents and assume that both of them start with service of the first among the waiting ones.  
Then  the trajectories of both the systems evolve in  exactly the same manner,  until the $\t$-queue gets empty.  There can be a change in the trajectories of the two systems,     upon a subsequent new $\t$-arrival.   We can have two scenarios as in Figure \ref{Figure_Threesystems}.  If $\i$-agents are absent at the  $\t$-arrival instance  in the original ${\cal O }$ system  (as in  sub-figure b), 
then again, both the systems continue to evolve in the same manner.  On the other hand, if     $\i$-agents are   deriving   service  (as in  sub-Figure a), the service of $\t$ agent is delayed in  the original ${\cal O }$ system till the end of the ongoing $\i$-busy period. While the service starts immediately in ${\cal M}_L$ system.   
Then the  trajectories in the two systems continue with the same difference,  until the end of the next  $\t$-idle period.  
At this point the difference:  a) either  gets   reduced, if the   $\t$ arrival marking the end of $\t$-idle period occurs after   sufficient time and  finds no $\i$-agent;  b) or   can increase, if the  $\t$-arrival occurs again during an $\i$-busy period;  c) or    can continue with almost previous value, if  the  $\t$-arrival occurs  immediately and finds  no $\i$-agent. And this continues. Thus   the  sojourn times  in ${\cal M}_L$ system are lower than or equal to that  in ${\cal O}$ system in all sample paths.   As we notice the difference between the two systems is because of   $\i$-busy cycles and this difference may diminish if the later  shorten. 
We will show this indeed is true in coming sections. 

\paragraph{${\cal M}_U$ system}
 Consider another $M/G/1$ system whose service times equal  $\Upsilon_\t + \Psi_\i$, where $\Psi_\i$ is an additional $\i$-busy period independent of $\Upsilon_\t$.
It is clear that this system dominates the ${\cal O}$ system everywhere (see ${\cal O}$ and  ${\cal M}_U$ trajectories in Figure  \ref{Figure_Threesystems}).  Hence the sojourn times of $\t$-agent in ${\cal O}$ system are upper bounded by that in ${\cal M}_U$ system (in all sample paths).
 Thus the expected sojourn time  of  ${\cal O} $ system is sandwiched 
  as below: 
  
  \vspace{-8mm}
{\Draft{}{\small}\begin{eqnarray}
   E^{{\cal M}_L}[S_\t]   \le E^{{\cal O}}[S_\t] \le  E^{{\cal M}_U}[S_\t] .
  \end{eqnarray}}

 \vspace{-4mm}
\subsubsection{Performance of ${\cal M}_L$ and ${\cal M}_U$ systems} 
In Lemma \ref{Lemma_Tt}, we obtained the first two moments of the $\i$-busy period and the EST, $\Upsilon_\t$. 
Using the well known formula for the expected sojourn time of an $M/G/1$ queue,  we have:

\vspace{-3mm}
 {\Draft{}{\small}$$E^{{\cal M}_L}[S_{\t}] = E[\Upsilon_{\t}] + \frac{\lambda_{\t} E[\Upsilon^{2}_{\t}]}{2(1-\rho^{{\cal M}_L}_{\t})} \mbox{ {\normalsize with}  } \rho^{{\cal M}_L}_{\t} = \lambda_{\t} E[\Upsilon_{\t}].$$}
Similarly with $\rho^{{\cal M}_U}_{\t} = \lambda_{\t} E[\Upsilon_{\t} + \Psi_\i]$,
{\Draft{}{\small}\begin{eqnarray*}
E^{{\cal M}_U}[S_{\t}] = E\big[\Upsilon_{\t} +\Psi_\i \big ]  
 +
 \frac{\lambda_{\t} \left ( E[\Upsilon^{2}_{\t}]   + E [\Psi_\i^2] + 2 E[\Psi_\i] E\Upsilon_\t]\right ) }{2(1-\rho^{{\cal M}_U}_{\t})}. 
 \end{eqnarray*}}
From Lemma \ref{Lemma_Tt} constants $\{c_i\}$, moments of busy period $E[\Psi_\i]$, $E[\Psi\up2_\i]$  converge to zero as $\mu_\i \to \infty$, and so the difference $E^{{\cal M}_U}[S_{\t}]  - E^{{\cal M}_L}[S_{\t}] $ converges to zero.    
In fact this is true even 
 when $\mu_\i, \lambda_\i$  jointly converge to  $ \infty$ while maintaining $\rho_\i  =\lambda_\i / \mu_\i$ constant.  If $\mu_\i \to \infty$ for a fixed $\lambda_\i$, then the load factor also decreases to zero in limit. Thus the result would have been true only for low load factors. But by maintaining the ratio $\rho_{\i}$ fixed when $\mu_\i \to \infty$, we ensured that {\it the approximation is good for any given  load factor and for any given admission control $p$, i.e., for any $(\rho_{\i}, p)$.} 
 Under SFJ limit, using Lemma  \ref{Lemma_Tt}:
 
 \vspace{-5mm}
 {\Draft{}{\small}\begin{eqnarray} \label{ES_p}
E_{PS}[S_\t(p)] &:= &E^{{\cal O}}_{PS}[S_\t(p)]  \approx     
\Draft{
\frac{ a_0} {\mu_\t } + \frac{  \lambda_\t a_0^2} { \mu_\t^2 \left ( 1- \frac{ a_0\lambda_\t} {\mu_\t } \right ) }
= 
}{}
  \frac{1}{\tilde{\mu}_{\t,p} (1- \tilde{\rho}_{\t,p})},  \\
\mbox{ with } \tilde{\rho}_{\t,p} &=& \rho_{\t} a_0,  \,\,\, \tilde{\mu}_{\t,p} = \frac{ {\mu_\t}}{  a_{0}} \mbox{ and } \rho_\t := \frac{\lambda_\t}{\mu_\t} . \nonumber
\end{eqnarray}}
 \Dtls{}{ In the above, where $\tilde{\rho}_{\t,p}$ is the inflated load factor. The above is similar to the sojourn time of an $M/M/1$ queue with equivalent service rate   
{\vspace{-4mm}$$  \tilde{\mu}_{\t,p} = \frac{ {\mu_\t}}{  a_{0}}=  \frac{ {\mu_\t}}{  a_{0} (p) }.$$}
  }
  Thus the achievable region under SFJ limit is given by:
  
   \vspace{1mm}
\hspace{-2mm}\fbox{\begin{minipage}{8.cm}\vspace{-3.mm}{ \small \begin{eqnarray*}
{\cal A}_{PS} = \bigg \{ \left ( (1-p)+ \frac{  p (\rho_{\i,p})^K }{ a_{0}}, \,\,   \frac{a_{0}}{\mu_\t  \left ( 1- a_{0} \rho_{\t} \right )}   \right )
\Draft{}{&\\ & \hspace{-32mm}}   \Big | \mbox{ with }    a_{0}\rho_\t  < 1, \  0 \le p \le 1   \bigg \}.
\end{eqnarray*}} \vspace{-4mm}\end{minipage} }\vspace{1mm}\\
In the above, condition  $a_{0} \rho_\t < 1$  ensures stability.

\subsection{Validation of Pseudo conservation law (\ref{Eqn_PC}),  Completeness}
By direct substitution\footnote{
By (\ref{Eqn_PBp}), $$1 - \rho_\i (1-P_B^{PS} ) = 1 / a_0$$ and so (see equation (\ref{Eqn_PC})) 
 $$\big (\mu_\t \big[1 - \rho_\i (1-P_B^{PS} ) \big ] - \lambda_\t \big )^{-1} \mbox{ equals} E_{PS}[S_\t],$$ given by  (\ref{ES_p}).\vspace{2mm}} one can verify that the performance measures of  $\beta^{PS}_{p,K}$ scheduler, for every $(p, K)$,  satisfy the pseudo conservation law (\ref{Eqn_PC}).  
 
 Further as $K $ increases to $\infty$, the   blocking probability  $P^{PS}_B (1)$,  given by  equation (\ref{Eqn_PBp}), decreases to zero if $\rho_\i \le 1$.
When $\rho_\i > 1$, using simple computations\footnote{It is easy to verify as $K \to \infty$ that:
\begin{eqnarray*}
\frac{\rho_\i^K}{\sum_{l=0}^K \rho_\i^l} = \frac{1}{\sum_{l=0}^K \rho_\i^{- (K-l)}} = \frac{1}{\sum_{l=0}^K \rho_\i^{-l}} \to
  \frac{1}{\frac{1}{1 - \rho_\i^{-1}}} = 1 - \frac{1}{\rho_\i}.
\end{eqnarray*}
}, one can show that  $$P^{PS}_B (1) \to 1 - 1/ \rho_\i$$ and only $p_B > 1 - 1/ \rho_\i$ can be a part of the ${\cal A}_{hetero}^{static}$.
 
  Also  it is easy to verify that  the function, $p \mapsto P^{PS}_B (p)$, is continuous in $p$ for any $K$. Thus by intermediate value theorem, all the points of ${\cal A}^{hetero}_{static}$ can be achieved by these schedulers. And hence the family of schedulers, 
 $${\cal F}^{PS}  :=  \bigg \{\beta^{PS}_{p,K}, \ 0 \le p \le 1, K  \bigg \}, \vspace{-1mm}$$  is   complete, when $\rho_\i \le 1$.  It is important to note here that these schedulers achieve the entire static region, nevertheless a larger $K$ implies a larger time spent by $\i$-agents in the system. Thus system may have a restriction on the size of $K$ to be used based on  other QoS requirements. 

\section{ Capacity Division ($CD$) Policies}
\label{sec_CD}

  In the previous    section, when an admitted $\i$-customer pre-empts the ongoing service of $\t$-customer, the entire system capacity is transferred to $\i$-customer. 
 In this \Thesis{section}{chapter} we analyze a different scheduling policy. Here the capacity is not completely transferred, but rather a fraction of it is used by each $\i$-customer. The $\t$-customer is  continued with the remaining capacity. 

 \subsubsection*{Service Discipline}  Each $\i$-customer uses {\small $(1/K)$}-th part of the capacity, $\mu_\i/K$. If the system has only one $\i$-customer, the remaining capacity i.e., {\small $(K-1)/{K}$}-th part of the capacity is utilized by the $\t$-customer. In other words, $\t$-class is served with rate {\small $\mu_\t (K-1) / K$}. If there are {\small $0 \le  l \le K$} number of $\i$-customers receiving the service, then {\small $({l}/{K})$}-th part 
 of the capacity is used by the $\i$-customers and the $\t$-customer is served at rate
 {\small $((K-l)/K) \mu_\t$}. This continues up to {\small $K$} $\i$-customers,  and any further (admitted) $\i$-arrival  departs without service. Whenever an existing $\i$-customer departs, the capacity is {\it readjusted to an appropriate higher value for $\t$-customer.}

It is more complicated to obtain the analysis of this model.  Now the effective server time depends upon the number of $\i$-customers in the system at the service start.  
However, the  analysis of $\i$-class is simpler for exactly the same reasons as in the previous model and this is considered first.

  \subsection*{Blocking Probability of $\i$-class}
\label{sec_Pb1}
It is clear that  $\i$-class evolution is again independent of $\t$-class evolution and its analysis is considered first.
    Consider any fixed $0 \le  p \le 1$. 
  The $\i$-inter arrival times  are exponentially distributed with parameter $\lambda_{\i}p$.
  Say there are $l$ $\i$-agents in the system  (note $l \le K$). Each one of them receive service at rate $\mu_\i/K$ and this happens simultaneously. Thus  the first  departure time would be exponentially distributed with parameter $l\mu_\i/K$.  
  This is again a continuous time Markov jump process and its transitions are as shown in 
  Figure  \ref{Figure_Transitions_1}.  In fact the $\i$-agents evolve like the well known, finite capacity and finite buffer queueing system, $M/M/K/K$ queue. The stationary distribution of such a queue  is well known and in particular  (see for e.g.,  \cite{Hoel}): 
\vspace{-1mm}
 \begin{eqnarray*}
\csn{\pi}_K &=&  \frac{(K\rho_{\i,p})^K}{ K! \csn{a}_0}  \:\: \mbox{where} \:\: \csn{a}_0 :=   \sum\limits_{j=0}^{K} \frac{\big (K \rho_{\i,p} \big )^j}{j!} .
\end{eqnarray*} 
As before, agents are admitted with probability 
  $(1-p)$, and hence the overall blocking probability by PASTA equals: \vspace{-1mm}
\begin{eqnarray} \label{PB_CD_P}
P_{B}^{CD}(p) &=& (1-p) + p \csn{\pi}_K 
= (1-p)  +  p\frac{( K \rho_{\i,p})^K}{K! \csn{a}_0} .
\end{eqnarray}
    \begin{figure}
    \vspace{5mm}
    \Draft{
\vspace{-19mm}
\hspace{-1.5mm}
\includegraphics[width=9cm,height=10cm]{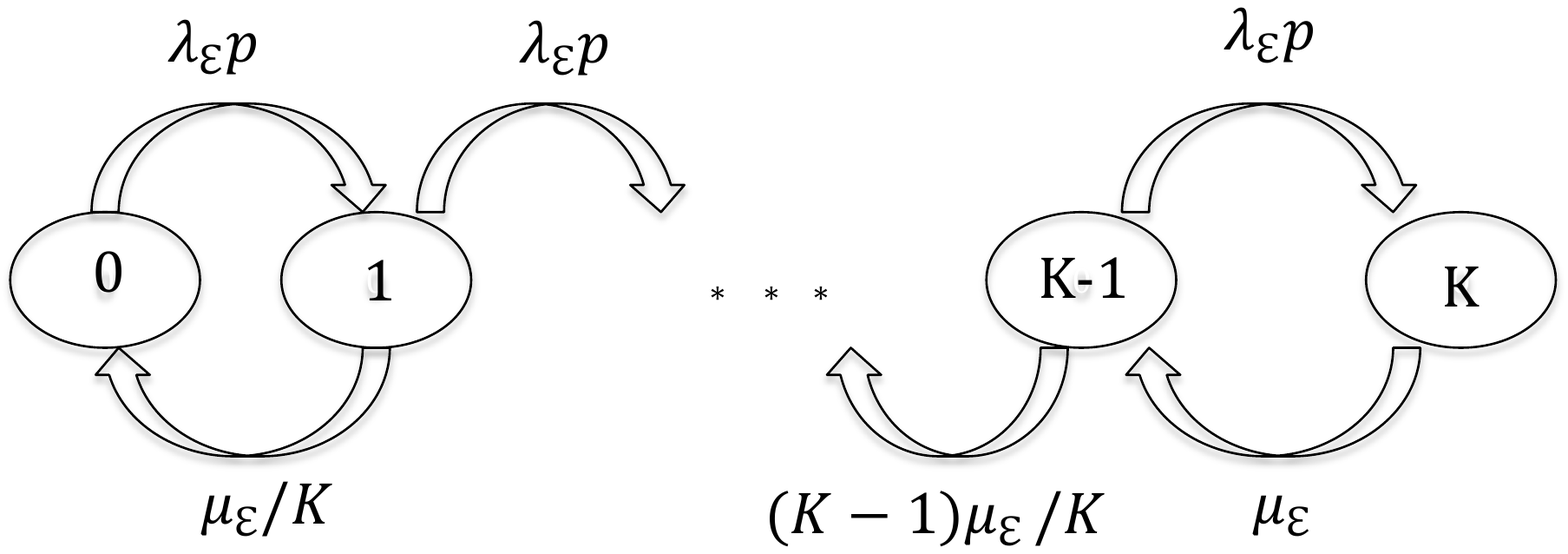}  
\vspace{-135mm}
\caption{ State transitions  for $\i$-agents in $CD$ model. \label{Figure_Transitions_1}}
\vspace{-5mm}}
{
\vspace{-12mm}
\hspace{-2.5mm}
\includegraphics[width=9cm,height=8cm]{NewCD_Transitions_exp} 
\vspace{-56mm}
\caption{ State transitions  for $\i$-agents in $CD$ model. \label{Figure_Transitions_1}}
\vspace{-5mm}}
\end{figure}

\subsection*{Expected sojourn time and Achievable region}
The expected sojourn time is obtained again using  dominating systems.  %
The idea  is once again to approximately  decouple the evolution of $\t$-agents from that of $\i$-agents. 
 The procedure  is similar,  however the current model is more complicated. Once again EST is denoted as $\csn{\Upsilon}_\t$,  has similar meaning as in \Thesis{section}{chapter} (\ref{sec_E_t}) and typical $\t$-sojourn time  equals the  sum of waiting time and the effective server time (EST), $\csn{(\Upsilon}_\t)$.

\subsubsection*{Analysis of effective server time $\csn{(\Upsilon}_\t)$} In the $CD$ model, EST is the total time period between the service start and service end of a typical $\t$-agent, when the capacity is divided possibly between the  two classes. 
With the arrival+admission of each $\i$-agent the server capacity available for the ongoing $\t$-agent reduces. With each $\i$-departure it increases, and  $\t$-service is completed in the midst of such rate changes. In fact, the  $\t$-agent's service is completely pre-empted with the admission of $K$-th $\i$-agent. The service would again be resumed, where left, when  one  of the $K$ $\i$-agents  depart. {\it The EST depends upon the number of $\i$-agents in the system at the service start.} Hence we introduce superscript $l$ in the notation of $\csn{\Upsilon}$. That is,  $\csn{\Upsilon}^l_\t$ represent the EST,  when it starts with $l $ $\i$-agents. 

Thus the analysis of EST for this model is not as easy as in $PS$ model. One can not estimate this using the number of interruptions and time per interruption as in $PS$ model. However, the underlying transitions are Markovian in nature, and hence we obtain the analysis by directly considering the EST's $\{ \csn{\Upsilon}_\t^l \}_l$. 
We have the following,   (proof in   Appendix B):

\vspace{2mm}
\begin{thm} 
\label{thm_Tt1}
 The first two moments of EST $\csn{\Upsilon}^0_\t$ ($l=0$) are:
 
\vspace{-3mm}
{\Draft{}{\small}\begin{eqnarray} \label{csnUpsilon1}
 E[\csn{\Upsilon}^0_\t]  
 &=&  \frac{ \csn{a}_0 + O (1/\mu_\i)} {\eta \mu_\tau + O (1/\mu_\i)}, \,\, \Draft{}{  \csn{a}_0 :=   \sum\limits_{j=0}^{K} \frac{\big ( \rho_{\i,p} \big )^j}{j!}} \mbox{ and }     \\ 
E \big [ \big (\csn{\Upsilon}^0_\t \big )^2 \big ] 
&
= &\frac{2 \frac{\csn{a}^2_0}{\eta \mu_\t} + O (1/\mu_\i)} {\eta \mu_\tau+ O (1/\mu_\i)} 
\mbox { with } \Draft{\csn{a}_0 :=   \sum\limits_{j=0}^{K} \frac{\big ( \rho_{\i,p} \big )^j}{j!} \mbox{ and } }{}  \eta  := \sum_{j=0}^{K-1} \frac{\big ( \rho_{\i,p} \big )^j}{j!}  \frac{K-j}{K}  , \nonumber
\end{eqnarray}}where $f (\mu_\i)  = O (1/\mu_\i)$ for any function $f$ implies, $f(\mu_\i)  \mu_\i \to  constant$ as  $\mu_\i \to \infty$, with $\rho_\i$ fixed.   {\hfill $\blacksquare$} 
\end{thm}

\begin{figure}
 \vspace{-8mm}
 \hspace{-6mm}
\includegraphics[width=10cm,height=11cm]{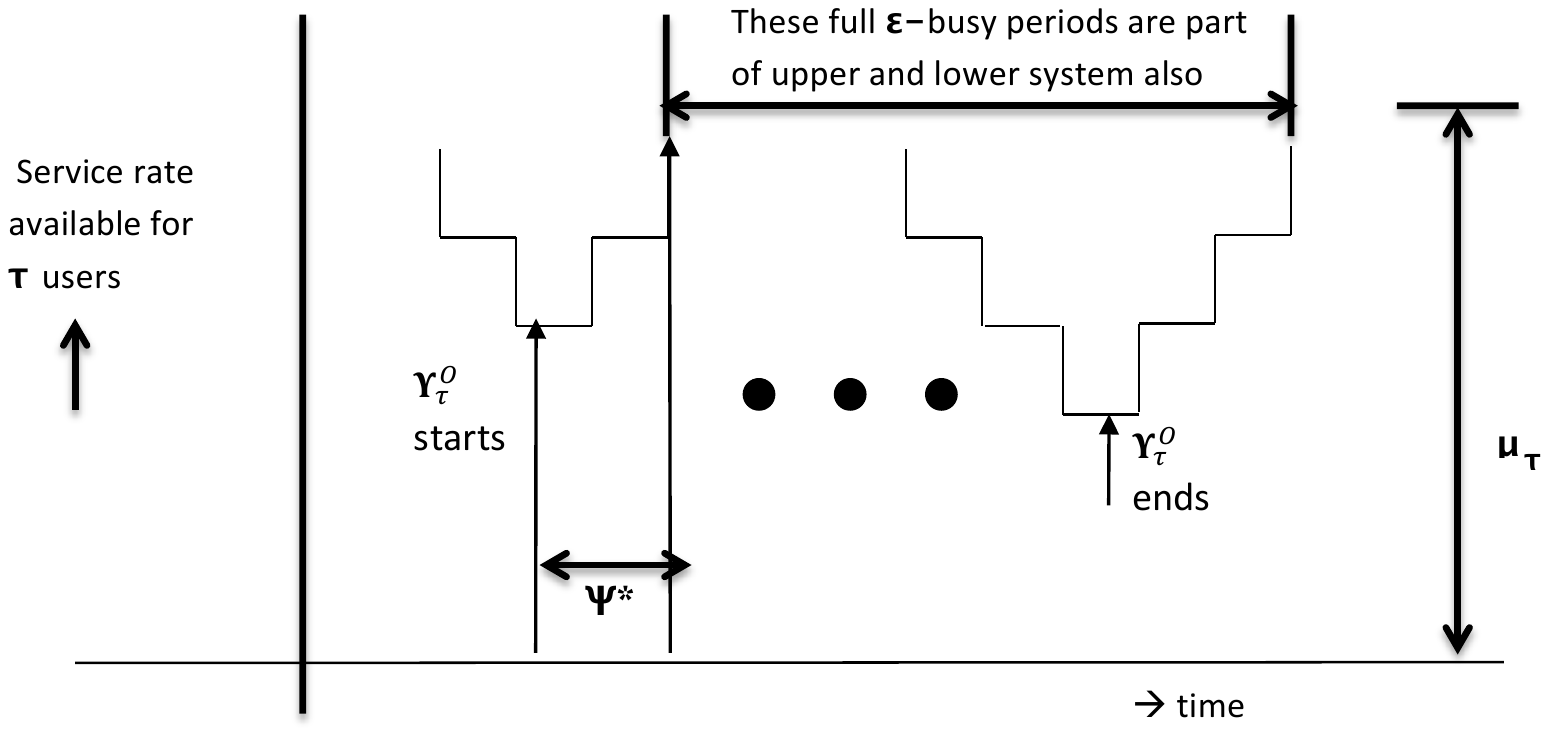}
\vspace{-75mm}
\caption{ Effective server time,  $\csn{\Upsilon}_\t$, in $CD$ Model) \label{Figure_uplowDominate} }
 \vspace{-5mm}
\end{figure}
 
 \vspace{2mm}
\paragraph*{Dominating systems}
It was not difficult to obtain the conditional moments of the EST, $\{ \csn{\Upsilon}_\t^l \}_l$. However to obtain the unconditional moments, \underline{one requires the stationary distribution of the $\i$-number $l$ at} \underline{service start of a typical $\t$-agent.} And this  again is not an easy task.  However  the various conditional moments differ from each other  (mostly) at maximum in one $\i$-busy period (see Figure \ref{Figure_uplowDominate}).  Hence one can possibly   obtain the (approximate) unconditional moments, along with M/G/1 queue approximation, using the idea of  dominating fictitious queues.

We again  construct two dominating systems, whose IID service times `dominate'  either side of the sequence $\{ \csn{\Upsilon}^{l_n}_{\t, n} \}_n$.   We first discuss the upper bounding system. The service times in original system $\csn{\Upsilon}^{l_n}_{\t, n} $ (for any $n$) can start and end in between  $\i$-busy period(s) as in Figure \ref{Figure_uplowDominate}.  Further these residual busy periods are correlated, for example the starting residual $\i$-busy period (call this as $\Psi^*$) is correlated with  $\csn{\Upsilon}^{l_{n-1}}_{\t, n-1} $  of previous $\t$-customer.  
To dominate any  $\csn{\Upsilon}_{\t}^{l}$  of original system with an IID version, we first   replace $\Psi^*$ with  busy period ${\tilde \Psi}^*$  of a CD system with $2K$ servers (each of capacity $\mu_\i/K$), when started with $K$ $\i$-customers and such that: a) if $l$ number of $\i$-customers are deriving service at the beginning of  $\Psi^*$, the residual service times (which are again exponential with   parameter $\mu_\i/K$, because of memoryless property) of those $l$ customers 
also equal the service time requirements of the first $l$ customers of
the $2K$ system;  b) the service times of the remaining $(K-l)$  $\i$-customers are independent copies of the  exponential random variable with the same parameter; c) further  inter arrival times and service times of all the new $\i$-customers coincide with that in the original system; and d) if a customer is not accepted in original system, we consider an independent service time for that customer.  With this construction, an $\i$-customer departure during
$\Psi^*$  of the original system definitely marks a departure in $2K$ system also,  any customer accepted in original system is also accepted in the $2K$  system. Thus the busy period ${\tilde \Psi}^*$ of the $2K$ system dominates the residual $\i$-busy period $\Psi^*$ at the start of the $\t$-customer service, irrespective of the number, $l$, of $\i$-customers existing in the original system at the start of $\Psi^*.$
In other words,  this time (say corresponding to $n$-th user) is independent of the quantities related to all other ($\ne n$) $\t$-customers  of  the original system and dominates $\Psi^*$ of the $n$-th customer  almost surely.

The   above constructed ${\tilde \Psi}^*$ of  the $2K$ system forms the beginning part of the $n$-th $\t$-customer service time in upper system with additional details:   a)  
the $\t$-customer in upper system is not served at the beginning for a duration equal to  ${\tilde \Psi}^*$; b) the service of $\t$-customer in upper system starts with full $\i$-busy periods,  and we  assume these  equal the full $\i$-busy periods of the original system that interrupted  the $n$-th $\t$-customer's service;  c)  if extra $\i$-busy periods are required to complete the $\t$-job we add independent copies of the $\i$-busy periods, but we do not couple the $\i$-busy periods that interrupt the $(n+1)$-th $\t$-customer. It is clear that $\t$-customer spends more time in upper system than in the original system.

A lower dominating system is obtained by using exactly the same construction, but  here the $\t$-customers are served with full capacity during ${\tilde \Psi}^*$, constructed using the $2K$ system.  Thus clearly the $\t$-customers spend less time (almost surely) in the lower system. And further the difference between the two dominating systems converges to zero because ${\tilde \Psi}^*$ (the busy period of CD system with $2K$ servers) also converges to zero as in proof of  Theorem \ref{thm_Tt1} provided in  Appendix B.

\begin{figure*} 
\vspace{42mm}
\hspace{-14.mm}
 \begin{minipage}{13cm}
\vspace{-85mm}

\includegraphics[width=11.cm,height=16cm]{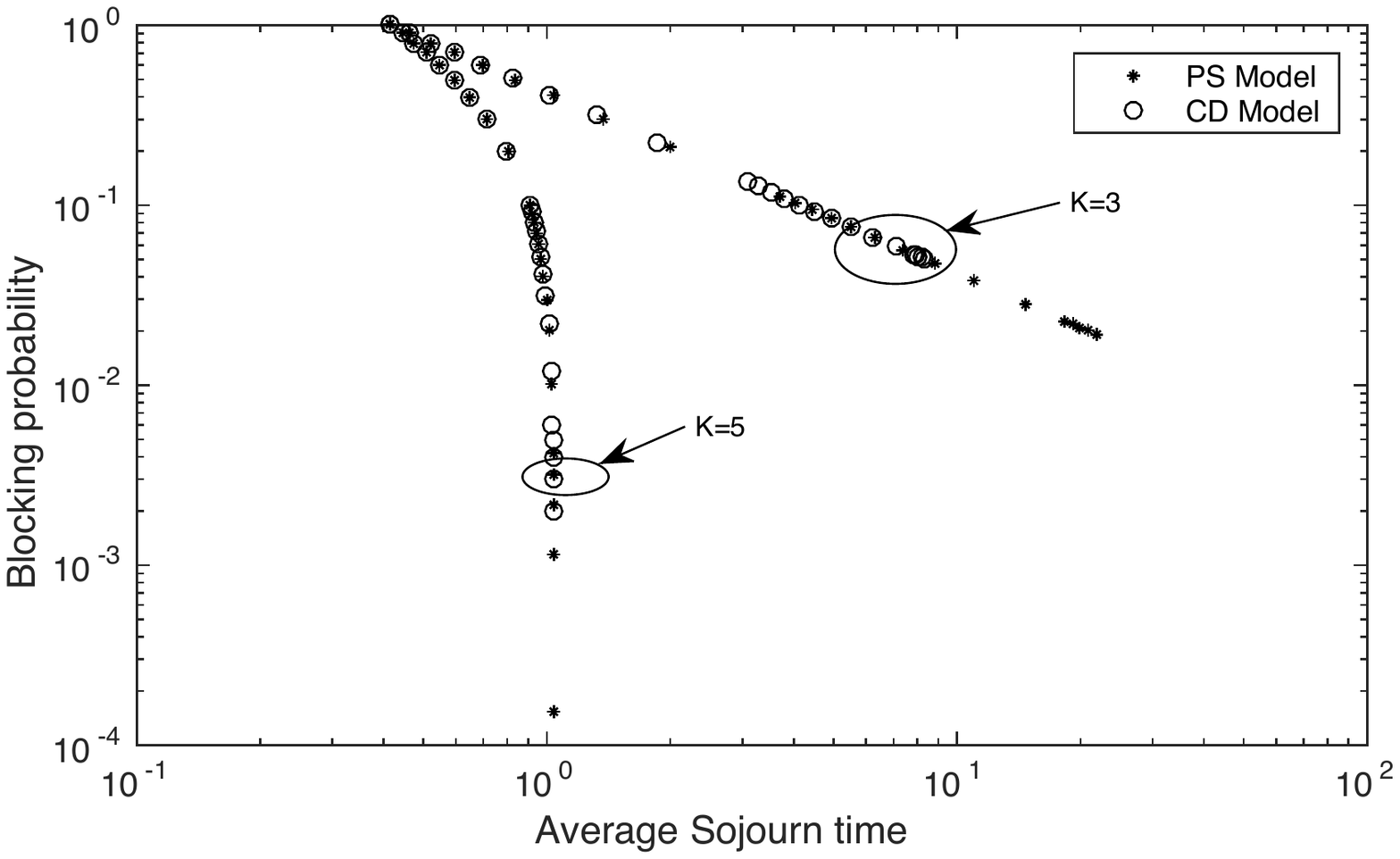}
\vspace{-54mm}
\caption{  Achievable regions ${\cal A}_{CD}$, ${\cal A}_{PS}$:  for  different $\rho_\i$, $\rho_\i = 0.9/K$  \label{Figure_ACS_PS}}
\vspace{20mm}
\end{minipage}
\hspace{-12mm}
\begin{minipage}{5cm}    
\vspace{-50mm}
\begin{minipage}{5cm} 
\vspace{-50mm}
\includegraphics[width=6.cm,height=6.cm]{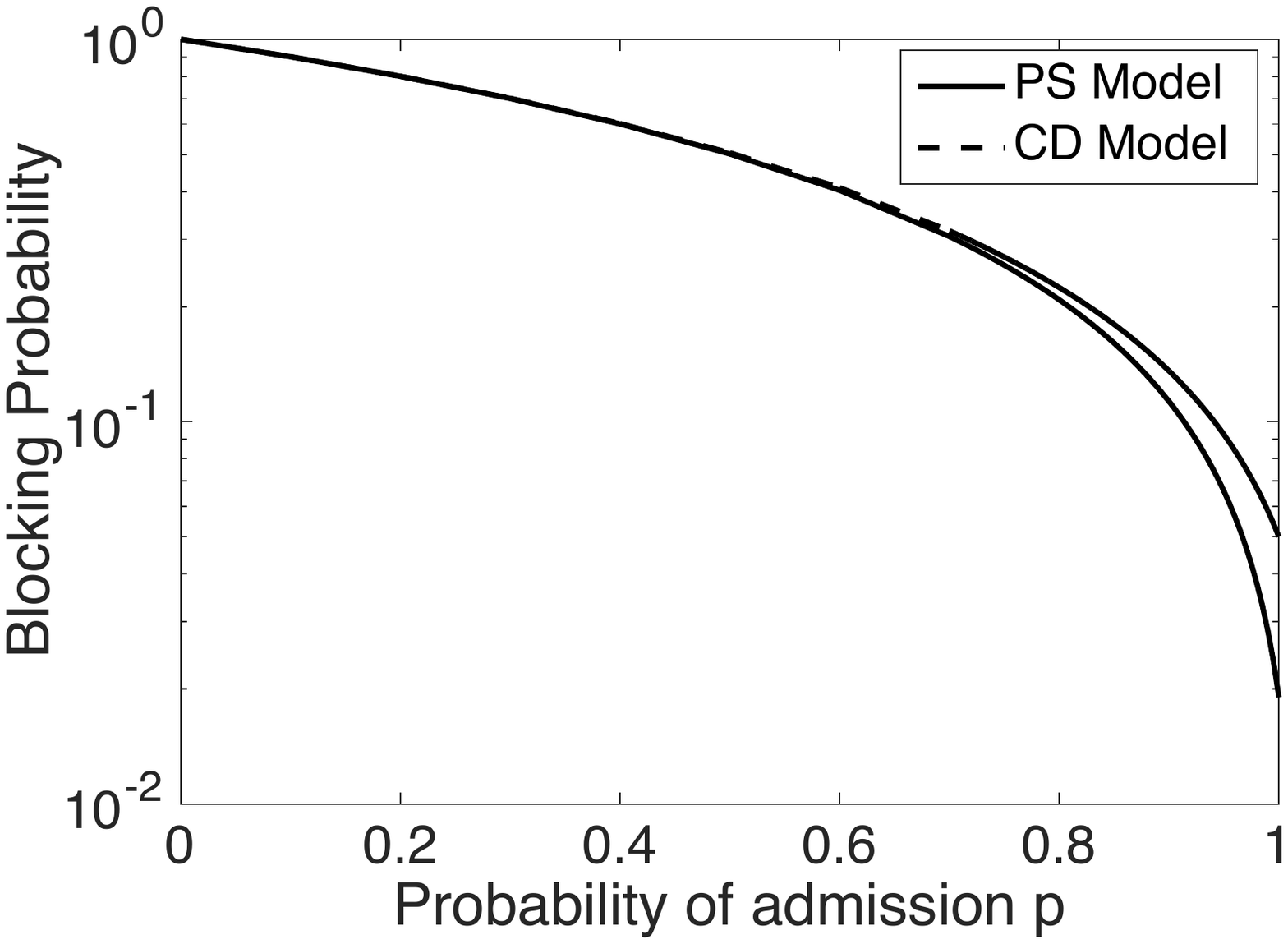}
\vspace{-25mm}
\caption{ Blocking Probability versus $p$ \label{Figure_PB_P}}
\end{minipage} 
\begin{minipage}{5cm}    
 \vspace{-17mm}
\includegraphics[width=6cm,height=6.cm]{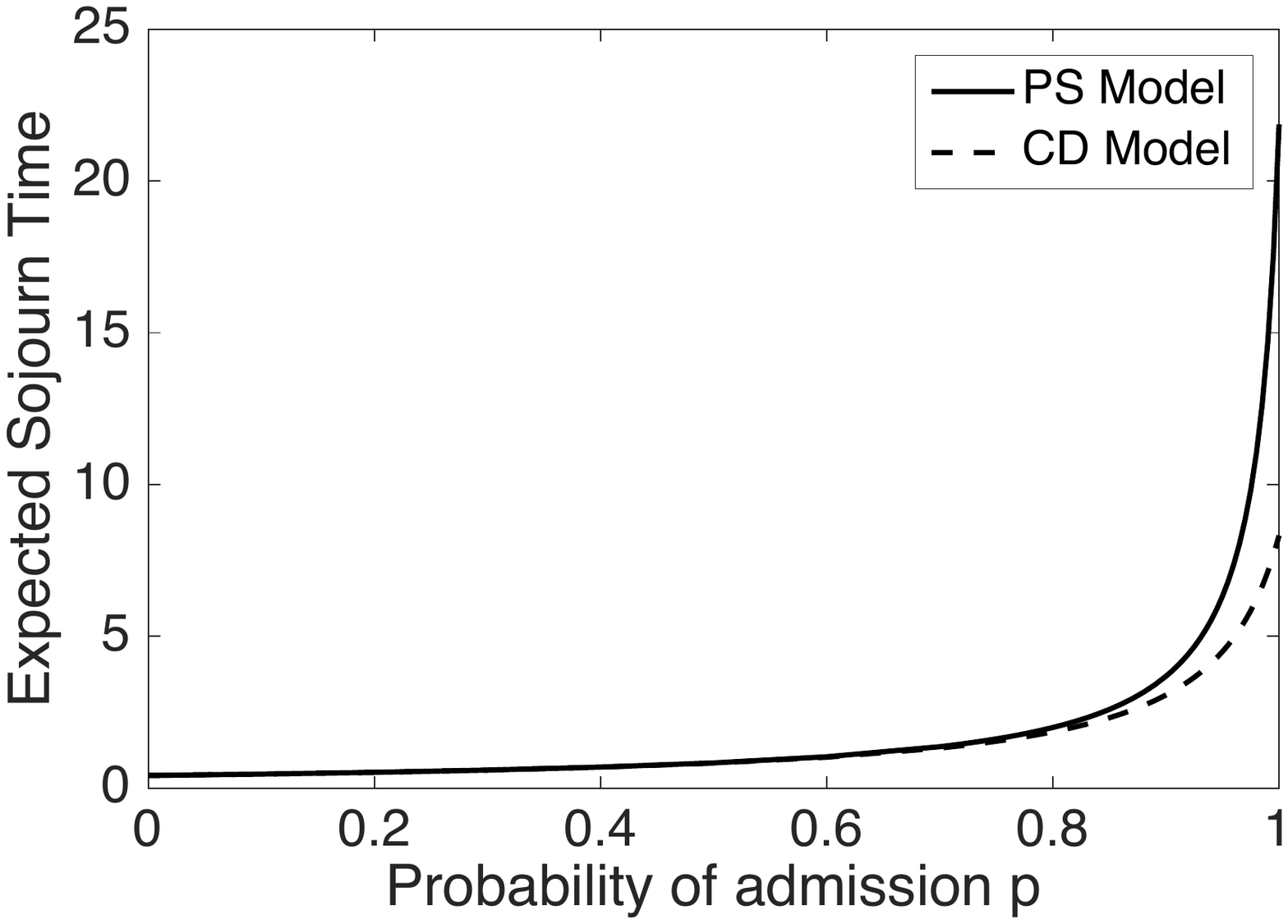}
\vspace{-23mm}
\caption{Expected Sojourn time versus $p$  \label{Figure_ES_P}}
 \end{minipage}   
\vspace{-25mm}
\end{minipage}   
  
\vspace{-16mm}
\end{figure*}

 \paragraph*{Performance}   Using exactly the same logic as in the previous model, one can show that the expected sojourn time of the $CD$ model can also be obtained as limit of  the expected sojourn times of   M/G/1 queues with service time moments given by that of $\csn{\Upsilon}^0_\t$ of Theorem \ref{thm_Tt1}.
To complete the analysis, one also needs to show the ergoidicity of the original system which is considered in \cite{Journ1}.
  Thus the 
  achievable region in the SFJ limit is given by   (Appendix B):
 
  \vspace{1mm}
  \hspace{-4mm}
\fbox{\begin{minipage}{8.7cm}\vspace{-3.mm}
{\small\begin{eqnarray*}
{\cal A}_{CD} = \Bigg \{ \left ( (1-p)  +  p\frac{(K \rho_{\i,p})^K}{K! \csn{a}_0}, \,\,   \frac{1}{\ddot{\mu}_{\t,p} \left(1 - \ddot{\rho}_{\t,p} \right)}   \right ) \Conf{}{& \\ &\hspace{-45mm} } :    \ddot{\rho}_{\t,p}  < 1, \:\: 0 \le p \le 1  \Bigg \}
 ,  \mbox{ with } \ddot{\rho}_{\t,p} = \frac{\lambda_{\t}}{\ddot{\mu}_{\t,p}},  \end{eqnarray*}\vspace{-3mm}
{\small $$
\csn{a}_0 :=   \sum\limits_{j=0}^{K} \frac{\big (K \rho_{\i,p} \big )^j}{j!}, \,\eta  := \sum_{j=0}^{K-1} \frac{\big (K \rho_{\i,p} \big )^j}{j!}  \frac{K-j}{K},\mbox{ and }  \ddot{\mu}_{\t,p} = \frac{\eta \mu_{\t}}{\csn{a}_{0}}.    $$}
}\vspace{-3.mm}\end{minipage}}  \vspace{1mm}

\NTR{

By direct substitution\footnote{From equation (\ref{PB_CD_P}), 
\vspace{-4mm}
$$1 - \rho_\i (1-P_B^{CD} ) =\frac{\sum_{j=0}^{K-1} \frac{\big (K \rho_{\i,p} \big )^j}{j!}  \frac{K-j}{K}}{  \sum\limits_{j=0}^{K} \frac{\big (K \rho_{\i,p} \big )^j}{j!}} = \frac{\eta}{\csn{a}_0} = \nu_K^{CD},$$ and so  $\big (\mu_\t  \nu_K^{CD}- \lambda_\t \big )^{-1}$ (see equation (\ref{Eqn_PC})) equals 
$$\Expp_{PS}[S_\t] =   \frac{1}{\ddot{\mu}_{\t,p} \left(1 - \ddot{\rho}_{\t,p} \right)} , $$ as given  in ${\cal A}_{CD}$.}  one can verify that  the $CD$ policies also satisfy the pseudo conservation law (\ref{Eqn_PC}). Further, this family 
 is also a complete family of schedulers for exactly the same reasons as used for PS policy and further using Lemma \ref{lemma_pb(1)} of Appendix B.

}{
By direct substitution,   we see that  the $CD$ policies also satisfy the pseudo conservation law (\ref{Eqn_PC}). They  also form a complete family of schedulers, for exactly the same reasons as that for $PS$ policy when $\rho_\i \le 1$ (details in \cite{Journ}).  For $\rho_\i > 1$, they cover only partial achievable region (\cite{Journ}).}

 \section{Numerical examples}
 \label{eqn_numerical}
 
\subsubsection*{Random system with large $\mu_{\i}$} We conduct Monte-Carlo simulations to estimate the performance of both the  policies. We basically generate random trajectories of the two arrival processes, job requirements and   study  the system evolution when it schedules agents according to $PS$/$CD$ policy.   We estimated the blocking probability and expected sojourn time for $\i$ and $\t$-agents respectively, using sample means, for different values of   $(p, K)$.

  In Figure \ref{Figure_Sim_thARegion}, we consider an example to compare the theoretical expressions with the ones estimated using Monte-Carlo simulations for $PS$ policy.  We consider two different values of $\rho_\i$. We  notice negligible difference between the theoretical and simulated values with $\mu_\i = 100$. However even with $\mu_\i = 20$, the difference is about 10-12$\%$ for most of the cases.  \NTR{}{Some more examples, including $CD$ policy, are available in \cite{Journ}.}
  
\subsubsection*{Achievable region} is also plotted in Figure \ref{Figure_Sim_thARegion} for different values of $\rho_\i$. 
Towards this, we plot $E^{PS}[S_\t(p)]$ versus $P_B^{PS}(p)$, for  $p \in \{ i \delta: 0 \le i \le 1/\delta \}$ with sufficiently small $\delta>0$.  It is a convex curve.  We notice a downward shift (improvement) in the curve with smaller $\rho_\i$, as anticipated.  However the formula derived, helps us understand the exact amount of shift.      We plotted the curves only in the $\t$-stability region, 
$\{ \lambda_\t : a_{0}\rho_\t  < 1\}. $  
{\it Unlike the case of homogeneous agents, the $\t$-stability region varies with the scheduling policy. This is because, varying fractions of $\i$-agents are lost  for different $p$, which can expand or contract the stability region.}

 \vspace{1mm}
	\subsubsection*{Comparison of the two  policies} 
We compare the achievable regions of $PS$ and $CD$ policies by plotting ${\cal A}_{CD}$ and ${\cal A}_{PS}$. We set $\rho_\i = 0.9 / K $, $\lambda_\t =5.6$ $ \mu_\t=8$ and $K = 3$ or $5$. 

In Figure \ref{Figure_ACS_PS}, we plot the achievable region for both the models/policies, i.e, we plot $E[S_{\t}(p)]$ versus $P_B(p)$, for different $p$. And in Figures \ref{Figure_PB_P} and \ref{Figure_ES_P}, we plot the performance measures $P_B (p)$ and $E[S_{\t}(p)]$ respectively versus $p$ with $K=3$. From Figure \ref{Figure_ACS_PS}, the two achievable (sub) regions overlap, however we observe from the Figures \ref{Figure_PB_P} and \ref{Figure_ES_P} that the performance measures of the two  models are different for the same $(p, K)$. But if we choose a $p$ and 	$p\prime$ such that $P_B^{CD}(p) = P_B^{PS} (p\prime)$, we observe that the two expected sojourn times are equal. Because of this the two achievable regions overlap in Figure \ref{Figure_ACS_PS}.   This  observation  {\it is precisely the pseudo-conservation law.  Whatever the policy used, once the blocking probabilities are the same the expected sojourn times are the same.} 

Now we will discuss a slightly different,  yet, a related important aspect. We would compare the two sets of policies, when $K$ (maximum number of parallel calls) is the same. 
As seen from the figures  the sub-achievable region of $CD$  policy,  with fixed $K$,  is a strict subset of that of the $PS$ policy. This is because the best possible blocking probability with $CD$ policy, 
 
 \vspace{-4mm}
{\small $$ \hspace{9mm}P_B^{CD}(1) =  \frac{( K\rho_{\i})^K / K!}{  \sum_{j=0}^{K}  {  (K \rho_{\i}  )^j}/{j!} }  
\ge  \frac{(\rho_{\i})^{K}}{\sum_{j=0}^K (\rho_\i )^j} = P_B^{PS} (1),
$$}is greater than that with the $PS$ policy. 
In Figure \ref{Figure_ACS_PS} the best $P_B$ with  $CD$ and $PS$ models/policies respectively is 0.002 and 0.0002 (0.05 and 0.019) when $K=5$ ($K=3$).
Thus it appears that the static achievable region would overlap for different policies, however the sub-regions covered by different policies can be different when $K$ is fixed.\Conf{ In future we would like to work on the entire possible achievable region, irrespective of the system configuration/model.}{}

 \begin{figure} [h]
 \vspace{2.5mm}
 \hspace{-2mm}
   \begin{minipage} {4.5cm}
\vspace{-15mm}
\includegraphics[width=5cm,height=7.cm]{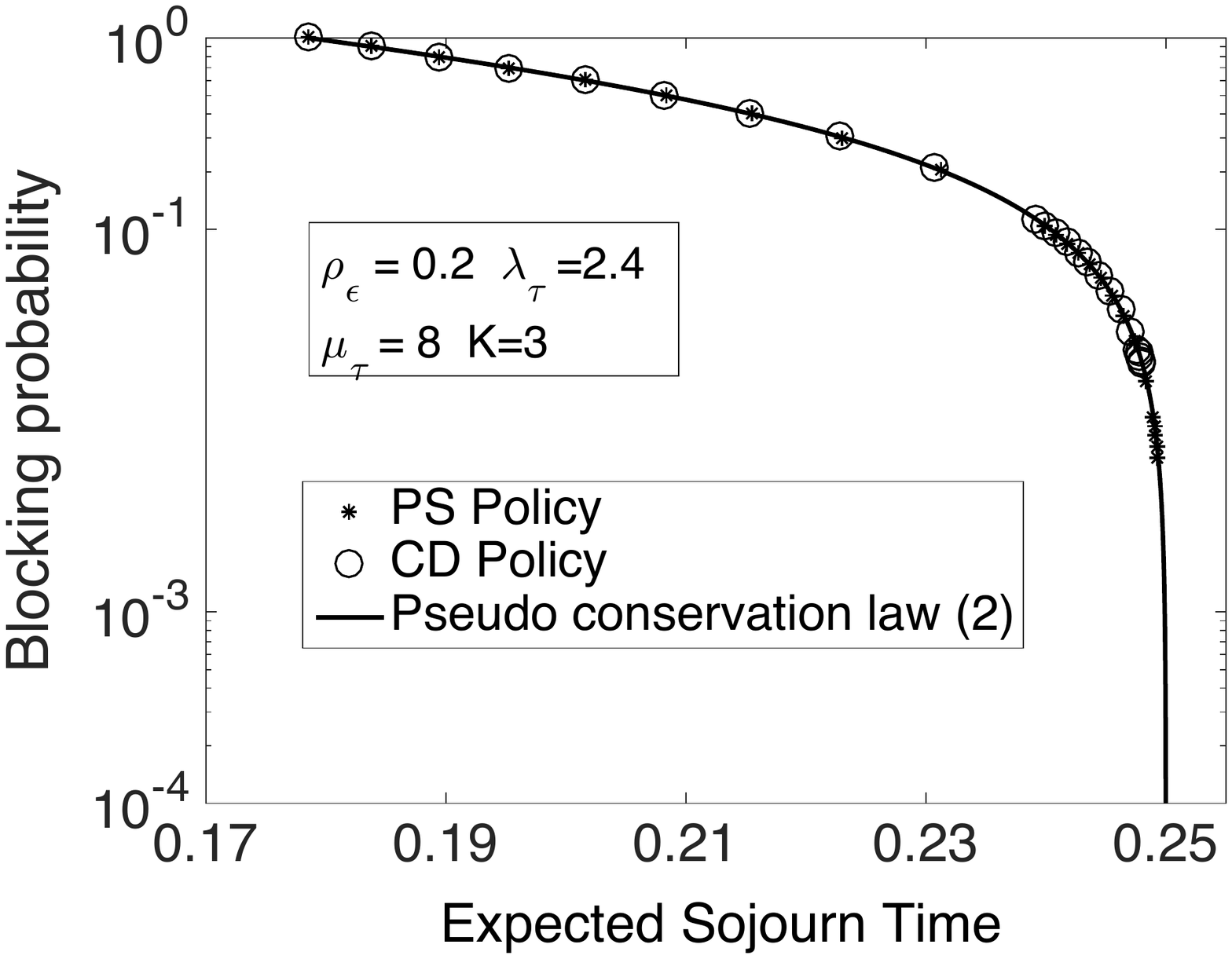}
\vspace{-31mm}
 \caption{  Static Achievable region  \label{Figure_EntireRegion}}
\end{minipage}
\hspace{1mm}
   \begin{minipage} {4.2cm}
\vspace{-28mm}
\includegraphics[width=4.8cm,height=9.cm]{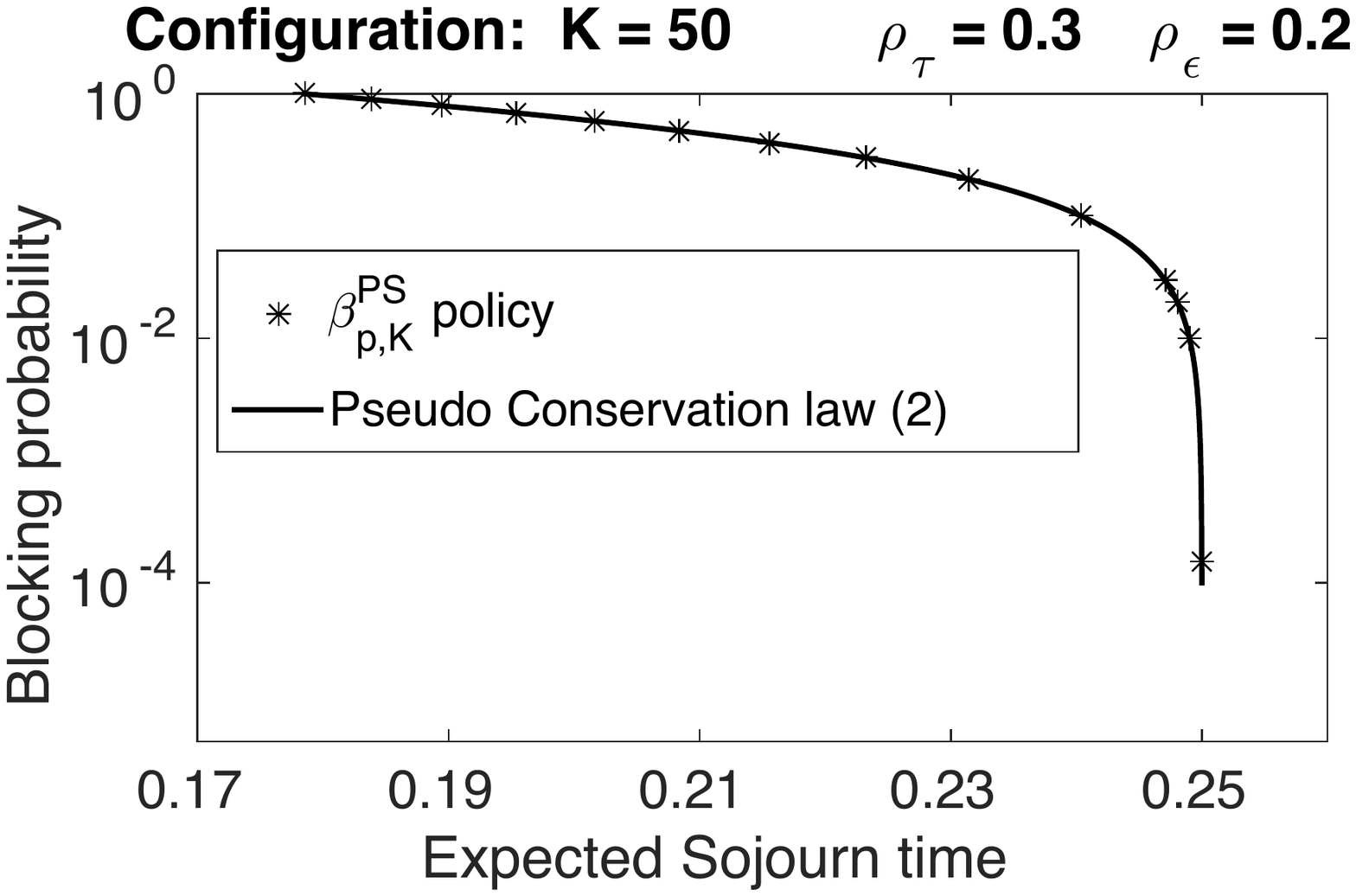}
\vspace{-39mm}
 \caption{  Completeness of  ${\cal F}^{PS}$    \label{Figure_PSEntireRegion}}
\end{minipage}
\vspace{-1mm}
\vspace{-1mm}
   \end{figure}

     \subsubsection*{Completeness} In Figure \ref{Figure_EntireRegion},  we plot pseudo-conservation law (\ref{Eqn_PC}). We also plot   the performance  of $PS$/$CD$ policies with $K= 3$ and for varying $p$.  We see that the three curves exactly overlap, again validating (\ref{Eqn_PC}).  
     For the same configuration we plot performance of $PS$ policies with a bigger $K = 50$,   in Figure \ref{Figure_PSEntireRegion}. With $K=50$ we are able to achieve a bigger part of the achievable region. One can achieve a similar result with $CD$ policy. With even bigger $K$ one can achieve further   lower parts of the pseudo-conservation curve.  However, as mentioned before, one may not be able to use a larger $K$ because of other QoS restrictions. For example, the $\i$ customers may not agree  for  a very small service rate ($\mu_\i/K$) which can prolong their stay in the system.  It is in this context that the $PS$ could be better than the $CD$ policies.  Even though both the sets of policies are complete, $PS$ policy achieves a bigger sub-region than the CD policy for the same $K$ (see Figures \ref{Figure_ACS_PS} and \ref{Figure_EntireRegion}).
 
\NTR{\section{A dynamic policy }
\label{sec_dynamic}

We consider dynamic policies (for $PS$ model) with an aim to demonstrate that the dynamic region is  bigger than the static region. Towards this we construct an example dynamic policy and show that the  block probability, for the same sojourn time $\Expp[S_\t]$, is better with the   dynamic policy.

The static policy of the previous sections   is modified  as follows. We refer this as  policy $\beta^d_p$. When there are no $\t$-agents in the system, i.e., during the $\t$-idle period, there is no admission control for $\i$-agents. An arriving $\i$-agent is admitted with probability one. Recall, however that service is offered to an admitted agent only when the number in system is less than   $K$.  When the system is in $\t$-busy period\footnote{Normally a busy period begins immediately with an arrival to an empty queue. However, in our system we say a $\t$-busy  period starts with the service start of that $\t$-agent, which arrives to an $\t$-empty queue.  
If $\i$-agents were present at the $\t$-arrival instance, the service of the $\t$-agent   is deferred till the end of the ongoing $\i$-busy period.}, i.e., when the $\t$-queue is non-empty, we admit the $\i$-agents  with probability $p$.   So, this is a dynamic policy which alternates between full  and partial admission.

Let $\Psi_\t$ and $\Id_\t$ respectively represent the busy and idle periods of the $\t$-agents. By stationarity, memoryless property, the consecutive busy, idle periods $\{\Psi_{\t, i} \}$,  $\{\Id_{\t,i} \}_i$  are independent and identically distributed. We have (proof is in Appendix C):
\begin{thm}
\label{thm_dynamic}
The block probability, $P^B_{d} (p)$,\:for the system with the dynamic policy $\beta^d_p$:

\vspace{-2mm}
{\Draft{\vspace{-4mm}}{\small}\begin{eqnarray}
\label{Eqn_PBpd}
\hspace{6mm} P^B_{d} (p) =  \frac{\Expp[\Id_{\t,1}] P^B (1)}{\Expp[\Psi_{\t,1}] + \Expp[\Id_{\t,1}]} +  \frac{\Expp[\Psi_{\t,1}] P^B (p) }{\Expp[\Psi_{\t,1}] + \Expp[\Id_{\t,1}]}   . \hspace{6mm}\blacksquare
\end{eqnarray}}
\end{thm}
  
Using the ideas of dominating systems as in the  section \ref{sec_static}  one can show that the moments of the idle,  busy periods of the original system with policy $\beta^d_p$ converges towards that of the equivalent $M/G/1$ system ${\cal M}_L$, as $\mu_\i \to \infty$.  Thus we will have for large values of 
$\mu_\i$:
{\small
\begin{eqnarray*} 
  \Expp[\Id_{\t,1}]  &\approx&   \frac{1}{\lambda_\t} ,    \\
\Expp^{\cal O}[\Psi_{\t,1}]  &\approx& \Expp^{{\cal M}_L}[\Psi_{\t,1}]  = \frac{   \Expp[\Upsilon_\t] }{1-  \lambda_\t \Expp[\Upsilon_\t]} 
 \to \hspace{1mm} \frac{a_0}{\mu_\t - \lambda_\t a_0} .
\end{eqnarray*}}
The second last equality is obtained using the  well known formula for the average busy period of an $M/G/1$ queue.
 It is easy to see that the sojourn time of the dynamic policy $\beta^d_p$ is same as that with static policy $\beta_p$ (asymptotically), while the blocking probability is improved from (\ref{Eqn_PBp}) to (\ref{Eqn_PBpd}). Note that $P^B(1) \le P^B (p)$ for any $p \le 1.$ Hence the dynamic policy performs better and the dynamic achievable region is bigger.  One can obtain similar improvement with $CD$ model.

 \begin{figure}
 \vspace{-1mm}
\begin{minipage} {9.6cm}  
\vspace{-6mm}
\includegraphics[width=8.5cm,height=9.cm]{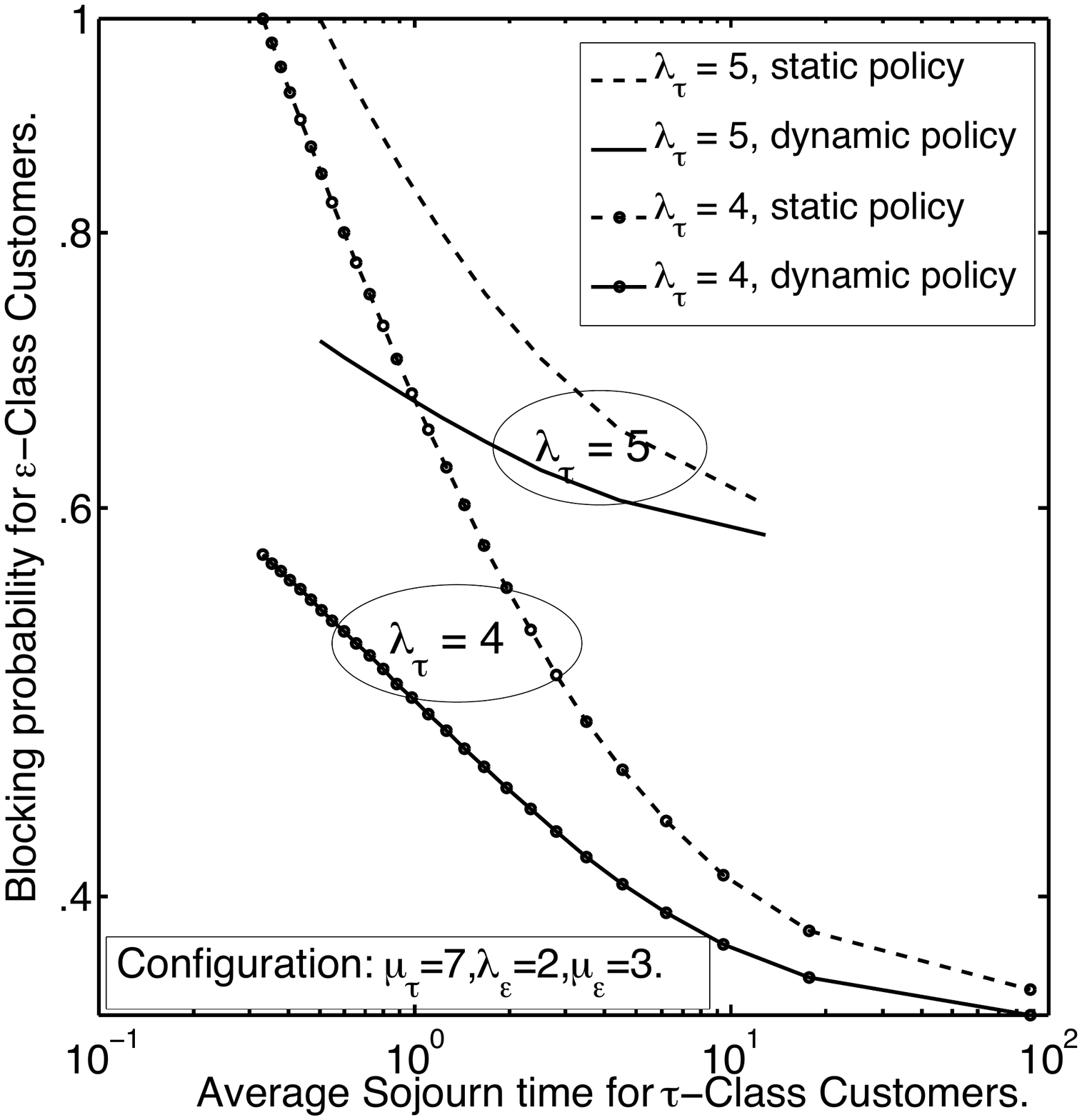}
\vspace{-15mm}
\caption{Static-dynamic policies, {\small$K = 4$}. \label{Figure_ARegion}} 
\end{minipage}
 \vspace{-2mm}
\end{figure}

 \subsection*{Numerical comparison of Dynamic and Static regions } 
In Figure \ref{Figure_ARegion}, we   compare   the performance   of the dynamic policy $\beta^d_p$ with the corresponding static policy, for $PS$ model. We notice a good improvement in the curve: blocking probability decreases significantly for the same expected sojourn time. This indicates that the dynamic region is strictly bigger than the static region, unlike the homogeneous case.  In future,
we would like to obtain   complete analysis of dynamic achievable region for this heterogeneous system.  }{} 
 
\vspace{-2mm}
\section{Conclusions and future work}
We consider a queueing system  with heterogeneous classes of agents. The impatient  class demands immediate service, hence receives the service immediately and if required in parallel with others. There is an admission control to ensure the QoS requirements of the other (tolerant) class.  The  tolerant class can wait for their turn, however would like to optimize their sojourn time.

We conjecture a pseudo conservation law for this lossy queueing system, which relates the blocking probability of impatient agents  to the expected sojourn time of the tolerant agents, in a short and frequent job (SFJ) limit-regime for the former.  The pseudo conservation law should be satisfied by all the  policies, that are  static (do not depend on $\t$-state)  and work conserving (left over server capacity is completely used when there is a customer)  with respect to the tolerant agents.

We consider two families of scheduling policies, which differ in the way the system capacity  is shared between the two classes. With processor sharing  policy the  entire system capacity is transferred to impatient customer, once admitted. In the second policy, which we refer as capacity division policy,  only a (fixed) fraction of capacity is transferred to each admitted impatient customer.   

We obtain closed form expressions for the asymptotic performance measures, under SFJ limit,  for both the families of policies.  The two families satisfy the pseudo-conservation law.  Further, both the families are complete, i.e., they  attain every point of the   achievable region given by the pseudo-conservation law.
The $CD$ achievable region is a strict subset of the $PS$ region, when   restricted  to the same number of parallel service possibilities. This demonstrates the limitation of $CD$ model, which could be a more practically used model.  The $PS$ model can attain a much smaller blocking probability.

The results are asymptotic and are accurate when the arrival-departure rates of the impatient class is large. Usually such customers have short frequent job requirements and hence this is an useful asymptotic result.  Further, we have an upper and lower bound for the sojourn time performance, even when the rates are not large. 


 Towards the end, we briefly discuss dynamic policies. These are the policies that depend upon the state of both the classes. We derive asymptotic performance of  an example dynamic policy and establish that the dynamic region is strictly bigger than the static region.

\section*{Appendix A: Proof of Lemma \ref{Lemma_Tt}}   
By conditioning on $B_\t$, one can verify that 
{\footnotesize \begin{eqnarray*}E [N(B_\t) ] &=& \frac{\lambda_\i p} { \mu_\t}, \, E[ B_\t N(B_\t) ] = \frac{2 \lambda_\i p} { \mu_\t^2}, \\
E[   ( N(B_\t) )^2 ] & = & \frac{\lambda_\i p} { \mu_\t} + \frac{2 (\lambda_\i p)^2}{\mu_\t^2} .
  \end{eqnarray*}}
%
By conditioning on $N(B_\t)$ we obtain the first moment: 
{\footnotesize \begin{eqnarray}
 E[\Upsilon_{\t}] \hspace{-2mm}&=& \hspace{-2mm}E[B_{\t}]  + E \left [\sum\limits_{i=1}^{N(B_{\t})} \Psi_{\i,i} \right ] \\
 & & \hspace{-16mm} = E[B_{\t}] + E\left [E\left [\left . \sum\limits_{i=1}^{N(B_{\t})} \Psi_{\i,i} \right |{N(B_{\t})} \right ] \right ]  \nonumber 
 = \frac{1}{\mu_{\t}} + \frac{\lambda_{\i}p E[\Psi_{\i}]}{\mu_{\t}} .  
\end{eqnarray} }
Note   that the   busy periods $\{ \Psi_{\i, i} \}_i$ are IID.
 From  (\ref{Eqn_Tt})  we have: 
{\small \begin{eqnarray}
E[\Upsilon_{\t}^2] 
 &=& E[B_{\t}^2] + 2E \big [ B_{\t}\Upsilon_{\t}^{e} \big ] + E \big [(\Upsilon_{\t}^{e})^{2} \big ].
\end{eqnarray} }
By first conditioning on ($B_{\t}$, $N(B_{\t})$) and then on $B_\t$:  \vspace{-1mm}
{\footnotesize \begin{eqnarray}  
E \big [ B_{\t} \Upsilon_{\t}^e \big ] 
= E\left [E \left [ \left . B_{\t}\sum\limits_{i=1}^{N(B_{\t})} \Psi_{\i,i} \right  | B_{\t},N(B_{\t}) \right ] \right ] \nonumber  \\ 
= \lambda_{\i}p E \big [\Psi_{\i}]E\big [B_{\t}B_{\t} \big ] = \frac{2\lambda_{\i}pE \big [ \Psi_{\i} \big ]}{\mu^{2}_{\t}}.
\end{eqnarray}}
Conditioning as before and because of independence: 
{\footnotesize\begin{eqnarray}
E \big [(\Upsilon_{\t}^{e})^2 \big] & = & E\left[E\left[ \left . \bigg ( \sum\limits_{i=1}^{N(B_{\t})} \Psi_{\i,i} \bigg )^{2} \right |N(B_{\t})\right] \right], \nonumber \\ 
\Conf{
& & \hspace{-19mm} =  E\Bigg[E \big [ \Psi^{2}_{\i} \big ] N(B_{\t}) + (E \big [ \Psi_{\i}\big ])^{2} N(B_{\t}) \big ( N(B_{\t}) - 1 \big )\Bigg] \nonumber \\
\TR{ \\
 & & \hspace{-19mm} =  \frac{\lambda_{\i}pE \big[ \Psi^{2}_{\i} \big ]}{\mu_{\t}} + \big ( E \big [\Psi_\i\big ]\big )^{2}E\bigg [ E \bigg [ ( N(B_\t ) )^2 - N (B_\t) |B_{\t}\bigg ] \bigg ] \nonumber }{}\\
 }{}
 &=&  \frac{\lambda_{\i}pE\big [ \Psi^{2}_{\i}\big ]}{\mu_{\t}} + \frac{2 (\lambda_\i p)^2}{\mu_\t^2}{\big ( E \big [ \Psi_{\i} \big ]\big )^{2}} \nonumber
\TR{ \\
& & \hspace{-19mm} =  \frac{\lambda_{\i}p}{\mu_{\t}} \left[E\big [ \Psi^{2}_{\i}\big ] + \frac{2 \lambda_\i p}{\mu_{\t}} \big (E\big [\Psi_{\i}\big ]\big )^{2} \right]   . \hspace{7mm}}{,}
\end{eqnarray}}which simplifies to (\ref{Eqn_EPsiUpsilon}).

\subsubsection*{Busy period of $\i$-class}
 Busy period of any class is defined as the time till the first epoch at which all the customers of that class have departed. Let  $\Psi_k$, represent the busy period of  $\i$-class, when it begins with $k$ number of customers. Note that $\Psi_\i =  \Psi_1$.
 In all the discussions below, an arrival is meant an admitted arrival.

The busy period  $\Psi_1$ starts with the arrival of one $\i$-customer.  If the customer leaves before the next arrival, the busy period ends. On the other hand, if an arrival occurs before the departure of the existing customer, it marks the beginning of a busy period with two customers,   $\Psi_2$.
As seen in \Thesis{section}{chapter} \ref{sec_Pb} (see Fig. \ref{Figure_Transitions}), a departure time is memoryless, i.e., exponential random variable with parameter $\mu_\i$ irrespective of the number of customer sharing the service.  Let $D$ represent the departure time. The inter arrival time, $A$, is exponential with parameter $\lambda_\i p$. 
Let $W : = \min \{ D, A\}$ represent the minimum of the two. 
With these definitions\TR{, one can write the following equation which relates  $\Psi_1$ and $\Psi_2$:}{:}

\vspace{-5mm}
{\small
\begin{eqnarray}
  \Psi_1 = {1}_{\{D<A\}}\:  0 + {1}_{\{A<D\}} \:  \Psi_2 + W. \label{Eqn_Psi1}
\end{eqnarray}}
The busy period  $\Psi_2$ starts with two $\i$-customers.  If one of the two customers leave before the next arrival, it marks the beginning of the busy period $\Psi_1$, and an early arrival  marks the beginning of a busy period with three customers,   $\Psi_3$.
From Lemma \ref{lemma_exp}  of \Thesis{section}{chapter} \ref{sec_Pb},  the departure time of the earliest customer among the two is again   exponential random variable with parameter $\mu_\i$. Thus this departure time is also distributed as $D$, defined above.  
 The inter arrival time $A$ obviously remains the same as in the previous paragraph. Thus,  
$
  \Psi_2 = {1}_{\{D<A\}}\:  \Psi_1 + {1}_{\{A<D\}} \:  \Psi_3 + W.
$
Continuing using similar logic we have:
 
 \vspace{-5mm}
 {\small\begin{eqnarray}
  \Psi_{K} &=& {1}_{\{D<A\}}\: \Psi_{K-1}  + {1}_{\{A<D\}} \:  \Psi_{K} + W \mbox{ and }  \label{Eqn_Psii} \\
  \Psi_{i} &=& {1}_{\{D<A\}}\: \Psi_{i-1}  + {1}_{\{A<D\}} \:  \Psi_{i+1} + W \, \,\, \forall \, 1< i < K .\nonumber 
\end{eqnarray}}
In the first equation of (\ref{Eqn_Psii}) the two $\Psi_K$ are different, independent of each other, but they are identically distributed. For ease of notation, we represent them by the same symbol. Note, in all, that the random variables $W$, $D$ and $A$ have same stochastic nature and are correlated. Further,  if an arrival occurs before departure 
 when the system already has $K$ customers, the arrival is dropped.  By memoryless property of exponential distributions, we again have busy period $\Psi_K$.
Taking expectation of equations (\ref{Eqn_Psi1}) - (\ref{Eqn_Psii}) and solving  backward recursively ($\{a_i\}$, $\{b_i\}$ given in (\ref{Eqn_constants})):  

\vspace{-4mm}
{\small\begin{eqnarray}
E[\Psi_{i}] = \frac{ ia_{i} + b_{i}}{\mu_{\i}}  \mbox{ for all }  1 \le i \le K.
\label{Eqn_EPsi}
\end{eqnarray}}
Squaring and taking the expectation of  (\ref{Eqn_Psi1})   we obtain:

\vspace{-4mm}
{\small \begin{eqnarray*}
E \big[\Psi_1^2 \big ]  &=&  c_1 +  q E \big [ \Psi_2^2 \big ] \mbox{ where } q := E \big [ A < D \big ] \\
c_1 &=& 2 E \left [ W 1_{\{A < D\} } \right ] E \big[\Psi_2 \big ] +  E\big [ W^2 \big ] \\ 
& =&  \frac{2\lambda_{\i}p}{(\lambda_\i p + \mu_\i )^2} E\big [ \Psi_2 \big ] + \frac{2}{ (\lambda_{\i}p + \mu_\i)^2}. 
\end{eqnarray*}}
Terms $c_1$, $q$ simplify as in (\ref{Eqn_constants}).
Similarly from (\ref{Eqn_Psii})  we have
{\small \begin{eqnarray*}
E[\Psi_i^2 ] \hspace{-2mm}  &=& \hspace{-2mm} c_i +  q E \big [ \Psi_{i+1}^2 \big ] + (1-q) E \big [\Psi_{i-1}^2 \big ]   \mbox{ with } \\
c_i \hspace{-2mm}&=&\hspace{-2mm} {\small 2 E \left [ W 1_{\{A < D\} } \right ] E [ \Psi_{i+1}] } + 2 E [ W 1_{\{D < A\} } ] E [ \Psi_{i-1}] \\ 
&&{\small \hspace{10mm} +  E[W^2] }  \hspace{3mm} \mbox{ for any } 2 \le i  < K. 
\end{eqnarray*}}
Constant $c_i$ simplifies as in (\ref{Eqn_constants}). Now  squaring $\Psi_K$ of  (\ref{Eqn_Psii}):
{\small\begin{eqnarray*}
E \big[\Psi_K^2 \big ] &=& E \big[\Psi_{K-1}^2 \big ] + \frac{c_{K}}{(1-q)}. 
\end{eqnarray*}}Solving the expressions backward recursively we obtain:
{\small \begin{eqnarray}
E[\Psi^{2}_\i] = E[\Psi_1^2] \hspace{-20mm} \\ &=& \frac{q^{K-1}\:c_{K}}{(1-q)^{K}} + \frac{q^{K-2}\:c_{K-1}}{(1-q)^{K-1}} + ... + \frac{q \:c_{2}}{(1-q)^{2}}  + \frac{c_{1}}{1-q}.  \hspace{4mm} \blacksquare \nonumber
\end{eqnarray}}

\vspace{-5mm}

 \NTR{

   \hspace{4mm}
   \vspace{40mm}

{\Large Appendix B and C  are  in the next pages.}

\newcommand{\Ki}[1]{[{#1}]}

\newcommand{\ItwoK}{[2, .., K-1] }

\newcommand{\II}[2]{{ I}_{#1}^{#2}}

\newcommand{\lamEpsp}  {{\tilde{\lambda}}} 

\newcommand{\Oe}{O_\i}

\newcommand{\Oesq}{O^{\mbox{\tiny (2)}}_{\i}}
 
 \onecolumn
\section*{Appendix B}

\noindent{\bf   Proof of Theorem \ref{thm_Tt1}:}

Let $\II{l}{m} :=[l, \cdots, m]$ represent the interval of integers, the big O notations are  shortly represented  by
{\small\begin{eqnarray}
\label{Eqn_notation}
\Oe := O(1/\mu_\i),  \,\, \Oesq :=O(1/\mu_\i\up2), \,\, \lamEpsp :=\lambda_\i p   \mbox{ and } \Ki{i} := K- i. \end{eqnarray}}  
 Let $\csn{\Upsilon}_l = \csn{\Upsilon}^l_\t$, with {\small$l \in \II{0}{ K-1}$}, represent a simpler notation for  EST when it begins with $l$  $\i$-agents. Let us begin with Êthe analysis of $\csn{\Upsilon}_0$.  If $\t$-agent leaves \Draft{the system }{}before  next  $\i$-arrival, the EST ends. If instead an $\i$-agent arrives  before, it marks the beginning of EST, \Draft{with one $\i$-agent, }{}$\csn{\Upsilon}_1$.Ê  Let $D_{\t}$ represent the departure time of $\t$-agent  and  $D_{\t} \sim \exp( \mu_\t)$. 
 It equals $B_\t$ since the service is offered at full capacity. 
 Let $A\sim \exp( \lamEpsp)$ represent the exponential inter arrival time  of admitted $\i$-agent. Let $\bar{W_0} : = \min \{ D_{\t}, A\}$ represent the minimum\Draft{ of the two}{}. Then,

\vspace{-4mm}
{\Draft{\vspace{-7mm}}{\small}
\begin{eqnarray}
  \csn{\Upsilon}_0 = {1}_{\{D_\t<A\}}\:  0 + {1}_{\{A<D_\t\}} \:  \csn{\Upsilon}_1 + \bar{W_{0}}. \label{Eqn_Upsilon0} 
\end{eqnarray}}

\vspace{-5mm}
     Let $D_{\i}\sim \exp( \mu_\i)$ represent the departure time of an $\i$-agent.  EST, $\csn{\Upsilon}_l$, starts with $l$ $\i$-agents. If one of the $\i$-agents depart before the next $\i$-arrival or $\t$-departure, it marks the beginning of EST $\csn{\Upsilon}_{l-1}$   and an early $\i$-arrival begins $\csn{\Upsilon}_{l+1}$. The $\t$-departure ends the EST. Let $D^l_{\i}$ represent the departure time of the earliest among  $l$ $\i$-agents, and note $D^l_{\i} \sim \exp( l\mu_\i)$. Let  $D^l_{\t}$ represent the departure time of $\t$-agent,  when  the capacity is shared with $l$ $\i$-agents, then    $D^l_{\t} \sim \exp(  \Ki{l}\mu_\t/K)$.  Inter arrival time, $A$, remains the same.
 Let $\bar{W_l} : = \min \{A , D^l_{\i}, D^l_{\t}\}$ represent the minimum of three, which is again exponentially distributed. Thus  for any $l \in \II{1}{K-1}$ we have:
 
 \vspace{-4mm}
{\Draft{\vspace{-6mm}}{\small}\begin{eqnarray}       \label{Eqn_Upsilonl}
\csn{\Upsilon}_l =\hspace{-1mm} {1}_{\{D^l_{\i} = \bar{W_l}\}}\:  \csn{\Upsilon}_{l-1} + {1}_{\{A= \bar{W_l}\}}\:  \csn{\Upsilon}_{l+1}  + \: \bar{W_l}.
\end{eqnarray}}
For the case with $K$ $\i$ agents, if one of them leave before next arrival, an EST with ($K-1$) $\i$-agents begins and an early arrival begins another EST with $K$ $\i$-agents.  Further, any arrival of $\i$-agent is dropped in this case. By memoryless property\Draft{ of exponential distributions,}{,} we again have busy period $\csn{\Upsilon}_K$ ($\bar{W_K }: = \min \{ D^K_{\i}, A\}$):

\vspace{-4mm}
{\Draft{}{\small}\begin{eqnarray}
  \csn{\Upsilon}_K = {1}_{\{D^K_{\i}<A\}}\: \csn{\Upsilon}_{K-1}  + {1}_{\{A<D^K_{\i}\}} \:  \csn{\Upsilon}_K + \bar{W_{K}}Ê
 . \label{Eqn_UpsilonK} 
\end{eqnarray}}
In the above the two $\csn{\Upsilon}_K$ are different, but  are identically distributed. For ease of notation, we represent them by the same symbol. 
Taking expectation of \Draft{equations }{}(\ref{Eqn_Upsilon0})-(\ref{Eqn_UpsilonK}) (\Draft{recall }{}$\Ki{i} := K-i$):

\vspace{-4mm}
{\Draft{}{\small }\begin{eqnarray}  \label{Eqn_zeroexp} 
E \big [ \csn{\Upsilon}_{0} \big ] &=& \frac{\lamEpsp}{\alpha_{0}} E \big [ \csn{\Upsilon}_{1} \big ]  + \frac{1}{\alpha_{0}},   \\
E \big [ \csn{\Upsilon}_{\Ki{i}} \big ]  &=&  \frac{ {\Ki{i}} \mu_\i}{\alpha_{\Ki{i}}} E \big [ \csn{\Upsilon}_{\Ki{i}-1} \big ] + \frac{\lamEpsp}{\alpha_{\Ki{i}}} E \big [ \csn{\Upsilon}_{\Ki{i}+1} \big ]  + \frac{1}{\alpha_{\Ki{i}}},\,\, \forall\,\, i \in \II{1}{K-1}  \nonumber \\
 E \big [ \csn{\Upsilon}_{K} \big ] &=&  \frac{K\mu_\i}{\alpha_K} E \big [ \csn{\Upsilon}_{K-1} \big ] + \frac{\lamEpsp}{\alpha_K} E \big [ \csn{\Upsilon}_{K} \big ]  + \frac{1}{\alpha_K}, \nonumber  
 \mbox{ where }  \\
 \alpha_{i}  &=&  \lamEpsp + i \mu_\i + \frac{\Ki{i} \mu_\t}{K} \mbox{ for any }   i \in \II{0}{K}
  .
 \nonumber
 \end{eqnarray}}
Solving the equations backward recursively (start with $K$\Draft{ and end till 1}{}): 	
\Draft{\begin{eqnarray} \label{Eqn_Ki}
E \big [ \csn{\Upsilon}_{K} \big ]  &=& m_K +E \big [ \csn{\Upsilon}_{K-1} \big ],  \nonumber  \\
E \big [ \csn{\Upsilon}_{K-1} \big ] &=&  m_{K-1} + n_{K-1} E \big [ \csn{\Upsilon}_{K-2} \big ] \mbox{ and }  \nonumber \\
E \big [ \csn{\Upsilon}_{\Ki{i}} \big ]  &=& m_{\Ki{i}} + n_{\Ki{i}} E \big [ \csn{\Upsilon}_{ \Ki{i}-1} \big ] \Draft{}{\nonumber}  \:\: \forall \:\:  i \in \II{2}{K-1}, 
\end{eqnarray}}
{{\small
\begin{eqnarray} \label{Eqn_Ki}
E \big [ \csn{\Upsilon}_{K} \big ] \Draft{=}{&=&} m_K +E \big [ \csn{\Upsilon}_{K-1} \big ],  \Draft{}{\nonumber  \\}
E \big [ \csn{\Upsilon}_{K-1} \big ] \Draft{=}{&=&} m_{K-1} + n_{K-1} E \big [ \csn{\Upsilon}_{K-2} \big ] \mbox{ and }  \Draft{}{\\}
E \big [ \csn{\Upsilon}_{\Ki{i}} \big ] \Draft{=}{&=&} m_{\Ki{i}} + n_{\Ki{i}} E \big [ \csn{\Upsilon}_{ \Ki{i}-1} \big ] \Draft{}{\nonumber}  \:\: \forall \:\:  i \in \II{2}{K-1}, 
\end{eqnarray}}}where the coefficients are defined recursively as below:
{\Draft{}{\small}\begin{eqnarray} \label{Eqn_Cns}
 m_K &=& \frac{1}{\gamma_k},  \,\, m_{K-1} = \frac{1+ m_K \lamEpsp}{\gamma_{K-1}}, \,\,  n_{K-1} = \frac{(K-1)\mu_\i}{\gamma_{K-1}},
  \,\,  \nonumber  \\
 m_{\Ki{i}}  &=& \frac{1+m_{\Ki{i}+1}\lamEpsp }{\alpha_{\Ki{i}} - n_{\Ki{i}+1} \lamEpsp}, \,\, \Draft{ \,\,\,}{}
 n_{\Ki{i}}  = \frac{(\Ki{i}) \mu_\i}{\alpha_{\Ki{i}} - n_{\Ki{i}+1} \lamEpsp}, \,\, \forall \:\: i \in \II{2}{K-1} \Draft{\mbox { where }}{} \nonumber \\
  \gamma_{i}  &=&  i \mu_\i + \frac{\Ki{i}\mu_\t}{K}, \:\: \forall \:\: i \in \II{0}{K}. 
\end{eqnarray}}
Using the first equation of (\ref{Eqn_zeroexp}) and  $\Expp[\csn{\Upsilon}_1]$ of equation (\ref{Eqn_Ki})  we obtain: 
 
{\Draft{}{\small}\begin{eqnarray} \label{Eqn_upsil} \hspace{5mm}
 \Expp[\csn{\Upsilon}_0] =   \frac{1 + m_{1} \lamEpsp}{\alpha_0 - n_{1} \lamEpsp}= \frac{1 + m_{1} \lamEpsp}{\lamEpsp + \mu_{\t} - n_{1} \lamEpsp}.
\end{eqnarray}}
   Squaring and taking the expectation of (\ref{Eqn_UpsilonK}) we get\footnote{Since the product of the two indicators is zero, we will not have cross correlation terms like $E \big [ \csn{\Upsilon}_{K1} \csn{\Upsilon}_{K2} \big ]$ etc. Furthe note that the indicators, $\{{\bar W}_l\}$ are independent of the ESTs $\{\csn{\Upsilon}_l\}_l$ on the right hand side of the equations 
   (\ref{Eqn_Upsilon0})-(\ref{Eqn_UpsilonK}). }
{\small \begin{eqnarray*}
E \big[ \big (\csn{\Upsilon}_K \big )^2 \big ] & =& \frac{\lamEpsp}{\alpha_{K}} E \big [\big ( \csn{\Upsilon}_{K} \big )^2  \big ]  + \frac{ K \mu_\i }{\alpha_{K}} E \big [\big ( \csn{\Upsilon}_{K-1} \big )^2  \big ]   \Draft{}{\\ && } + \frac{ 2 }{\alpha^2_{K}} + \frac{2 \lamEpsp}{\alpha^2_{K}} E \big [  \csn{\Upsilon}_{K}   \big ]  + \frac{ 2 K \mu_\i }{\alpha^2_{K}} E \big [ \csn{\Upsilon}_{K-1}   \big ]. 
\end{eqnarray*}}
Simplifying we obtain:
{\Draft{}{\small}
\begin{eqnarray}
E \big [ \big (\csn{\Upsilon}_K \big )^2 \big ]  &=& r_K + E \big [ \big (\csn{\Upsilon}_{K-1} \big )^2 \big ] \mbox{ with } \delta_K = \gamma_K = K \mu_\i\: \mbox{and} \\ 
r_K &=& \frac{2}{\delta_K \alpha_K} +\frac{ \sigma_K} {\delta_K }, \Draft{\\}{\,\,}
Ê\sigma_K \Draft{&:= &}{=}\frac{2 \delta_K}{\alpha_K} E \big [ \csn{\Upsilon}_{K-1} \big ] + \frac{2 \lamEpsp}{  \alpha_K} E \big [ \csn{\Upsilon}_{K} \big ]. \Draft{}{\nonumber} \label{Eqn_rk}
\end{eqnarray}}	
{Similarly from (\ref{Eqn_Upsilonl}) we obtain\Draft{:
   \begin{eqnarray*}
E \big [ \big (\csn{\Upsilon}_{K-1} \big )^2 \big ] &=& r_{K-1} + \frac{(K-1) \mu_\i}{\delta_{K-1}} E \big [ \big (\csn{\Upsilon}_{K-2} \big )^2 \big ], \mbox{where }  \nonumber \\ 
 r_{K-1} &=&  \frac{1}{\delta_{K-1}} \bigg [ \frac{2}{\alpha_{K-1}} + r_ {K} \lamEpsp + \sigma_{K-1} \bigg ], \delta_{K-1} = \gamma_{K-1},  \nonumber \\
 \sigma_{K-1}   &=&  \frac{2 \lamEpsp }{\alpha_{K-1}} \Expp[\csn{\Upsilon}_{K}] \: + \frac{2 (K-1) \mu_\i}{\alpha_{K-1}} \Expp[\csn{\Upsilon}_{K-2}]. \Draft{}{\nonumber \\ }
 \end{eqnarray*} 
Continuing in a similar way and using  (\ref{Eqn_Upsilonl}), for any  $i \in \II{2}{K-1},$}
{ for any  $i \in \II{1}{K-1},$}
 { \Draft{}{\small} \begin{eqnarray}
 E \big [ \big (\csn{\Upsilon}_{\Ki{i}} \big )^2 \big ] &=&  r_{\Ki{i}} + \frac{\Ki{i} \mu_\i}{\delta_{\Ki{i}}} E \big [ \big (\csn{\Upsilon}_{\Ki{i}-1} \big )^2 \big ] \mbox{ with}   \label{Eqn_Rec} \\  
 r_{\Ki{i}} &=&\frac{1}{\delta_{\Ki{i}}} \bigg [ \frac{2}{\alpha_{\Ki{i}}} + \lamEpsp r_{\Ki{i}+1} + \sigma_{\Ki{i}} \bigg ] , \nonumber \\
  \delta_{\Ki{i}} &=& \frac{\alpha_{\Ki{i}} \delta _{\Ki{i}+1} - \lamEpsp (\Ki{i}+1) \mu_\i}{\delta _{\Ki{i}+1}} \Draft{\mbox{\,\,and } \\}{,\nonumber \\ }
 \sigma_{\Ki{i}}  &=& \frac{2 \lamEpsp }{\alpha_{\Ki{i}}} \Expp[\csn{\Upsilon}_{\Ki{i}+1}] \: + \frac{2 (\Ki{i}) \mu_\i}{\alpha_{\Ki{i}}} \Expp[\csn{\Upsilon}_{\Ki{i}-1}].   \nonumber
\end{eqnarray}} }
{From (\ref{Eqn_Upsilon0}), \vspace{-8mm}
{\small $$ \hspace{15mm}
E \big [ \big (\csn{\Upsilon}_0 \big )^2 \big ]  =      \frac{\lamEpsp }{\alpha_0}  E \big [ \big (\csn{\Upsilon}_{1} \big )^2 \big ] +  \frac{2  }{\alpha_0\up{2}} +
\frac{2 \lamEpsp \Expp[\csn{\Upsilon}_1]}{\alpha_0\up{2}}.
 $$}
Further using equation
}
(\ref{Eqn_Rec}) with $i =K-1$ or $\Ki{i}=1$: 
{\Draft{}{\small} \begin{eqnarray}  E \big [ \big (\csn{\Upsilon}_0 \big )^2 \big ] = \frac{1}{\delta_0} \bigg [ \frac{2}{\alpha_0} + \lamEpsp r_1 + \frac{2 \lamEpsp \Expp[\csn{\Upsilon}_1]}{\alpha_0} \bigg ],\,\, \delta_{0} = \frac{\alpha_{0} \delta _{1} - \lamEpsp  \mu_\i}{\delta _{1}}
   .  \label{Eqn_meansqr}
   \end{eqnarray}}

 \noindent{\bf SFT Limit:}  From (\ref{Eqn_upsil}) and using  (\ref{Eqn_mn}) of Lemma \ref{Constant_n} (see (\ref{Eqn_notation})): 
{\Draft{}{\small}\begin{eqnarray} \label{Sigma}
 E \big [ \csn{\Upsilon}_0  \big ]  = \frac{ \csn{a}_0 + \Oe} {\eta \mu_\t + \Oe} .
\end{eqnarray}}
 Thus and considering the limit (forward) recursively in (\ref{Eqn_zeroexp})
{\Draft{}{\small} $$
  E \big [ \csn{\Upsilon}_l  \big ] = \lim E \big [ \csn{\Upsilon}_0  \big ] + \Oe = \frac{\csn{a}_0}{\eta \mu_\t} + \Oe.
 $$}
 Hence for all $ i \in \II{0}{K-1} $ from (\ref{Eqn_Rec})
{\Draft{}{\small} $$
  \sigma_{\Ki{i}} = \theta + \Oe, \,\,\,\,\, r_K \lamEpsp =   \frac{\rho_{\i,p} \theta}{K} + \Oe \,\, \mbox{where}\,\, \theta := \frac{2 \csn{a}_{0}}{\eta \mu_\t}.$$}
  Again considering limits  (backward) recursively in (\ref{Eqn_Rec}),  while using the above two equations  and equation (\ref{Eqn_deltai})  of Lemma \ref{Constant_n} and backward induction (as in Lemma \ref{Constant_n}) we obtain:
{\Draft{}{\small} \begin{eqnarray*}
  \lamEpsp r_{\Ki{i}}  &=& \left ( \frac{ \rho_{\i,p}  } { \Ki{i}} + \Oe \right )  \left (  \lamEpsp r_{\Ki{i}+1}   + \theta + \Oe  \right ) \\
 &=&  \left ( \frac{ \rho_{\i,p}  } { \Ki{i}} + \Oe  \right ) 
   \left ( \sum_{j=0}^{i-1} \frac{ \rho_{\i,p}^{j+1} \theta}{ (\Ki{i}+1) \cdots (\Ki{i}+1+j) }   
   + \theta + \Oe  \right) \\
   &=&   \sum_{j=0}^{i} \frac{ \rho_{\i,p}^{j+1} \theta}{ \Ki{i} \cdots (\Ki{i}+j) }   
   +   \Oe   \mbox{ for all } i \in \II{0}{K-1}.
\end{eqnarray*}}
Simplifying we obtain:

\vspace{-3.5mm}
{\Draft{}{\small}
 $$ \hspace{5mm} \lamEpsp r_1 + \frac{2 \lamEpsp \Expp[\csn{\Upsilon}_1]}{\alpha_0} = \frac{2 \csn{a}^2_{0}}{\eta\mu_{\t}} + \Oe  .  $$}
 Further using   $\delta_0$ of  (\ref{Eqn_deltai}) of Lemma \ref{Constant_n}  we obtain the asymptotic limit of the second moment (\ref{Eqn_meansqr}).  
 	 \hfill{$\blacksquare$}
 
  \begin{lemma} \label{Constant_n}
We have the following asymptotic results for the coefficients defined in the proof of Theorem \ref{thm_Tt1}:
{\Draft{}{\small}\begin{eqnarray} 
\hspace{-1mm} n_{\Ki{i}} &=& 1 - \frac{\mu_\t \omega_{\Ki{i}}}{\mu_\i} + \Oesq,\,\, \, i \in \II{1}{K-1} \mbox{ with}\, \label{Eqn_nki} \\
 \omega_{\Ki{i}} &:=& \frac{1}{K}\sum\limits_{j=0}^{i-1}\frac{(i-j)\rho_{\i,p}^{j}}{(\Ki{i}+j) (\Ki{i}+j-1) \cdots (\Ki{i})}, \nonumber\\
 \lamEpsp + \mu_{\t} - n_{1} \lamEpsp  &= &  \eta + \Oe , \hspace{5mm } 1 + m_1 \lamEpsp = \csn{a}_{0} + \Oe, \,\, \label{Eqn_mn} 
 \end{eqnarray}}
\Draft{\vspace{-7mm}}{\vspace{-4mm}}
 {\Draft{}{\small}\begin{eqnarray} 
   \delta_{\Ki{i}}  =   \Ki{i}\mu_\i + \Ki{i} \omega_{\Ki{i}} \mu_{\t} + \Oe  \hspace{1mm} \forall \: i \in \II{1}{K-1}
\mbox{ and } \delta_{0}    =  \eta \mu_\t +  \Oe. 
      \label{Eqn_deltai}   
 \end{eqnarray}}
\end{lemma}
{\bf Proof:} We begin with terms $\{n_{\Ki{i}}\}$ and prove the required result by backward mathematical induction.
 From  (\ref{Eqn_Cns}),

{\Draft{}{\small}$$n_{K-1} = 1- \frac{\mu_{\t}\omega_{K-1}}{\mu_{\i}} + \Oesq,\: \mbox{ with } \omega_{K-1} := \frac{1}{K(K-1)}.$$}
Assume the statement holds for  $i = l-1$,  i.e., say:
{\Draft{}{\small}\begin{eqnarray} \label{cons_n_{l-1}}
n_{K-l+1} &=& 1 - \frac{\mu_\t \omega_{K-l+1}}{\mu_\i} + \Oesq \mbox{ and \,} \Draft{}{\nonumber \\} 
\omega_{K-l+1} \Draft{:=}{&:=&} \frac{1}{K}\sum\limits_{j=0}^{l-2}\frac{(l-1-j)\rho_{\i,p}^{j}}{(K-l+1+j)\Draft{(K-l+1+j-1)}{} \cdots (K-l+1) }. \nonumber
\end{eqnarray}}
We need to prove the result for $i = l.$ From (\ref{Eqn_Cns}) and substituting the above
{\Draft{}{\small} 
\begin{eqnarray}
n_{K-l} &=& \frac{(K-l) \mu_\i}{\alpha_{K-l} - n_{K-l+1} \lamEpsp}  = 
\frac{(K-l)\mu_\i}{ (K-l)\mu_\i + \frac{l \mu_{\t}}{K} +  \mu_{\t} \rho_{\i,p} \omega_{K-l+1} + \Oe } \nonumber \\ \Draft{ \nonumber \\}{}
\Draft{&=& 1 - \frac{\frac{l \mu_{\t}}{K} + \mu_{\t} \rho_{\i,p} \omega_{K-l+1} +  \Oe }{(K-l) \mu_{\i} + \frac{l \mu_{\t}}{K} + \mu_{\t} \rho_{\i,p} \omega_{K-l+1} + \Oe  } + \frac{\frac{l \mu_{\t}}{K} + \mu_{\t} \rho_{\i,p} \omega_{K-l+1}}{(K-l) \mu_{\i}} - \frac{\frac{l \mu_{\t}}{K} + \mu_{\t} \rho_{\i,p} \omega_{K-l+1}}{(K-l) \mu_{\i}} 
\nonumber \\
\nonumber \\}{}
\Draft{&=& 1 - \frac{\frac{l \mu_{\t}}{K} + \mu_{\t} \rho_{\i,p} \omega_{K-l+1}}{(K-l) \mu_{\i}} +  \frac{(\frac{l \mu_{\t}}{K} + \mu_{\t} \rho_{\i,p} \omega_{K-l+1})^2 - (K-l)\mu_{\i}\Oe  + \Oe }{(K-l)\mu_{\i} \bigg( (K-l) \mu_{\i} + \frac{l \mu_{\t}}{K} + \mu_{\t} \rho_{\i,p} \omega_{K-l+1} +\Oe  \bigg )} \nonumber \\}{}
&=& 1 - \frac{\mu_\t \omega_{K-l}}{\mu_\i} + \Oesq \mbox{ as} \:\:\mu_\i \rightarrow  \infty \mbox{ with }   \rho_{\i}\:\: \mbox{constant.}                          \nonumber
\end{eqnarray}}
\underline{This proves  (\ref{Eqn_nki}).}
 It is easy to see that
{\Draft{}{\small} $$\omega_1 = \sum_{j=0}^{K-2} \frac{K-j}{K}  \frac{\rho_{\i,p}^j }{ j! }    \mbox{ {\normalsize and hence that} } \omega_1 \rho_{\i,p} + 1 = \eta , $$}where $\eta$ is defined in the hypothesis of the Theorem \ref{thm_Tt1}. Using this we obtain the \underline{first part of (\ref{Eqn_mn}):}

\vspace{-4mm}
{\Draft{}{\small}\begin{eqnarray}
\lamEpsp + \mu_{\t} - n_{1} \lamEpsp = (\omega_1 \rho_{\i,p} + 1) \mu_\t +  \Oe  =  \eta \mu_\t +  \Oe . \nonumber
\end{eqnarray}}
Using (\ref{Eqn_nki}),
 for all $i \in \II{2}{K-1}$ (note $\lambda_\i =\rho_\i \mu_\i$ with $\rho_\i$ fixed),
{\Draft{}{\small} \begin{eqnarray}
\alpha_{\Ki{i}} - n_{\Ki{i}+1} \lamEpsp &  = & \Ki{i} \mu_\i + \left  (\omega_{\Ki{i}+1} \rho_{\i,p} +  \frac{i}{K} \right )\mu_\t +  \Oe  \end{eqnarray}}  
\vspace{-2mm}
and hence \vspace{-5mm}
{\Draft{}{\small} \begin{eqnarray*}
\frac{ \lamEpsp}{
\alpha_{\Ki{i}} - n_{\Ki{i}+1} \lamEpsp }   =   \frac{ \rho_{\i,p} }{\Ki{i} }   +  \Oe . \nonumber
\end{eqnarray*}}
From  above and using a similar backward induction on $\{m_i\}$ of  (\ref{Eqn_Cns}) we  obtain, 

\vspace{-3mm}
{\Draft{}{\small}\begin{eqnarray} \label{Eqn_mi}
m_{\Ki{i}} \lamEpsp = \sum\limits_{j=1}^{i+1}\frac{\rho_{\i,p}^j}{(\Ki{i})(\Ki{i}+1) \cdots (\Ki{i}+j)}+  \Oe , \forall \:\: i \in \II{2}{K-1}. \:\:\end{eqnarray}}
Thus we get \underline{the second part of (\ref{Eqn_mn}):}  \vspace{-2mm}
{\Draft{}{\small}$$ \:\: 1 + m_1 \lamEpsp =  \sum\limits_{j=0}^{K} \frac{\rho^j_{\i,p}}{j!} = \csn{a}_{0} + \Oe .$$}
Let $\zeta_{\Ki{i}} := \Ki{i} \mu_\i / \delta_{\Ki{i}},$
 for all $i \in \II{1}{K-1}$. Using the recursive definition of $\delta_i$, as given  in (\ref{Eqn_Rec}), $\zeta_{\Ki{i}}$ satisfies the following recursive equation,
\vspace{-3mm}
{\Draft{}{\small}\begin{eqnarray}
 \zeta_{\Ki{i}} &=& \frac{\Ki{i} \mu_\i}{\alpha_{\Ki{i}} - \lamEpsp \zeta_{\Ki{i}+1}},
\end{eqnarray}}just like the recursive definition of  $\{n_i\}$ given in  (\ref{Eqn_Cns}). 
Further,   
 
{\Draft{}{\small} 
\begin{eqnarray*}
 \zeta_{K-1} &=&  (K-1)\mu_\i / \gamma_{K-1} =  n_{K-1}  \mbox{ {\normalsize  and hence using (\ref{Eqn_nki})} }  \\
 \zeta_{\Ki{i}} &=& n_{\Ki{i}} =  1 - \frac{\mu_\t \omega_{\Ki{i}}}{\mu_\i} + \Oesq  \mbox{ for all } i \in \II{1}{K-1}.
\end{eqnarray*}}
Thus we have the \underline{first part of (\ref{Eqn_deltai})}:
\begin{eqnarray*}
 \delta_{\Ki{i}} =  \Ki{i} \mu_\i +  \mu_\t  \Ki{i} \omega_{\Ki{i}}  + \Oe  \mbox{ for any } i \in \II{1}{K-1}.
\end{eqnarray*}
 And from (\ref{Eqn_meansqr}), \vspace{-7mm}
 $$\Draft{\hspace{60mm}}{\hspace{25mm}} \delta_{0}  = \mu_{\t} + \mu_{\t} \rho_{\i,p} \omega_{1} + \Oe   = \eta \mu_\t +  \Oe  . \Draft{\hspace{58mm}}{\hspace{8mm}} \blacksquare$$

\medskip

\begin{lemma} \label{lemma_pb(1)}
For any $\rho_\i > 1$,
$$
P_B^{CD} (1) \to  1 - \frac{1}{\rho_\i} \mbox{ and }   P_B^{CD} (1) \to 0 \mbox{ if } \rho_\i \le 1. $$
\end{lemma}
{\bf Proof:}  When $\rho_\i \le 1$ we have:
\begin{eqnarray*} P_{B}(1) = \frac{(K \rho_{\i})^K}{K! \sum\limits_{j=0}^{K} \frac{\big (K \rho_{\i} \big )^j}{j!}} = \frac{\frac{1}{\sum\limits_{j=0}^{K} \frac{\big (K \rho_{\i} \big )^j}{j!} }} {\frac{(K \rho_{\i})^K}{K!}}  &  =& \frac{1}{\sum\limits_{j=0}^K \frac{K (K-1) \cdots (K-j+1)}{K^{j}} \rho_{\i}^{j-K} } \\
&\leq&  \frac{1}{\sum\limits_{j=0}^K \frac{K (K-1) \cdots (K-j+1)}{K^{j}} }. \end{eqnarray*}
Let \begin{eqnarray*} f(K) &:=& \sum\limits_{j=0}^K \frac{K (K-1) \cdots (K-j+1)}{K^{j}}  \\
&=& 1+ (1-\frac{1}{K}) + (1-\frac{1}{K}) (1-\frac{2}{K}) + \cdots +\Pi_{i=1}^K (1-\frac{i}{K}) , \end{eqnarray*}
 and note that 
$$
f(K) \ge N \Pi_{i=1}^N (1-\frac{i}{K})   \mbox{ for any }N \le K.$$
Fix any $\varepsilon > 0.$
For any 
 $ N  $  there exist a large  $ K_N$ such that for all $K \ge K_N$,
  $$((1-\frac{1}{K+1})\cdots (1-\frac{N}{K+1}) \geq  (1-\varepsilon) \; \mbox{which implies}\; f(K) \geq N(1-\varepsilon), $$
 and hence $ f(K) \to \infty \;\mbox{as} \;K \to \infty. $
 Thus $ P_{B}^{CD}(1) \to 0. $
When $\rho_\i > 1$, it is clear after redefining  that:
$$
f(K) :=  \sum\limits_{j=0}^K \frac{K (K-1) \cdots (K-j+1)}{K^{j}} \rho_{\i}^{j-K}  \le\sum\limits_{j=0}^K  \rho_{\i}^{j-K} = \sum\limits_{j=0}^K  \rho_{\i}^{-j} $$
 hence 
$$
\lim_{K \to\infty}  f(K) \le \frac{1}{1 - \rho_\i^{-1}}.
$$ 
On the other hand for any $\varepsilon > 0$,  as before for any $N$ for all $K > K_N$ we have: 
$$
f(K) \ge \sum_{j=K-N}^K  \rho_{\i}^{j-K} (1-\varepsilon) = \sum_{j=0}^N  \rho_{\i}^{-j} (1-\varepsilon) .
$$ By first letting $N \to \infty$ we have $$
\lim_{K \to\infty}  f(K) \ge   (1-\varepsilon) \left ( \frac{1}{1 - \rho_\i^{-1}}  \right ) \mbox{ and then with } \varepsilon \to 0 
\lim_{K \to\infty}  f(K) \ge     \frac{1}{1-\rho_\i ^{-1}}   . \hspace{10mm}\Box$$

 }{}

 \section*{Appendix C:  Proof of Theorem \ref{thm_dynamic}}  
  For any static policy $\beta_p$,  the processor sharing system with   $\i$-agents    is ergodic. With   $L_p(T)$ representing the number of $\i$ agents lost in  time $T$, when the arrivals are admitted at rate $p$,  we have:

{\small
$$ \hspace{15mm}
\lim_{T \to \infty}  \frac{L_p(T)}{T}  =   P^B (p) \mbox{ \normalsize almost surely,}
$$}where $P^B (p)$ is given by equation  (\ref{Eqn_PBp}).

Let $L_p^d (T)$  represent the number of $\i$-agents lost in  time $T$ with dynamic policy.  Let $\Id_\t (T)$, $\Psi_\t (T)$ respectively represent the total $\t$-idle time and total $\t$-busy period until time $T$. These are basically the sum of all the busy/idle periods that elapsed till the time $T$. Note that  $\Id_\t (T) + \Psi_\t (T) = T.$ In this case, 

\vspace{-1mm}
{\small $$ 
 L_p^d (T) =  L_p^d ( \Id_\t (T) ) +  L_p^d ( \Psi_\t (T)).
$$}
At the end epoch of any busy period   ($\Psi_{\t,i}$  for some $i$), the system is completely empty. That is, agents of both the classes are absent. 
Also a $\t$-busy period starts only once the system is free from all of its  $\i$-agents.
 Thus there are no $\i$-agents in the system at both start and end epochs of a $\t$-busy period.  Hence the evolution of the  $\i$-class loss counting process that occurred during disjoint time intervals (of $\t$-busy periods) constituting $\Psi_\t (T)$ is stochastically
equivalent to the $\i$-class loss counting process  that would have evolved in a continuous time interval of length exactly $\Psi_\t (T)$.  
This is because of the memoryless property associated with Poisson arrival process.
Thus we have:

{\small $$ \hspace{9mm}
\frac{L^d_p ( \Psi_\t (T))}{\Psi_\t (T)} \to P^B (p) \mbox{ almost surely }. 
$$}
In a similar way at the start/end epoch of any $\t$-idle period the system is free of $\i$-agents.  
Using similar arguments, and  because all $\i$-agents are admitted during idle periods we have:

{\small $$
\frac{L^d_p ( \Psi_\t (T))}{\Psi_\t (T)} \to P^B (1) \mbox{ almost surely }. 
$$}
Further using renewal reward theorem, one can show that the following happens  almost surely:

{\small\begin{eqnarray*}
\frac{ \Psi_\t (T))}{T}  \to  \frac{\Expp[\Psi_{\t,1}]}{\Expp[\Psi_{\t,1}] + \Expp[\Id_{\t,1}]},  \mbox{  }  
\frac{ \Id_\t (T))}{T}  \to  \frac{\Expp[\Id_{\t,1}]}{\Expp[\Psi_{\t,1}] + \Expp[\Id_{\t,1}]}.  
\end{eqnarray*}}
Using all the results established so far,  equation (\ref{Eqn_PBpd}) follows because: 
{\small \begin{eqnarray*}
\hspace{25mm} P^B_d(p) = \lim_{T\to \infty} \frac{L^d_p(T)}{T}. 
\hspace{26mm}\blacksquare
\end{eqnarray*}}

\end{document}